\documentclass[11pt]{article}
\usepackage{lscape}
 \usepackage[ansinew]{inputenc}
 \usepackage[dvips]{graphicx}
\usepackage[psamsfonts]{amssymb} \usepackage[all]{xy}  \usepackage[all]{xy}
 \usepackage[psamsfonts]{eucal}
 \usepackage[psamsfonts]{eucal}
\usepackage{amssymb, amsmath}
\usepackage{amsthm}
\usepackage{geometry}
\usepackage{enumerate}
\usepackage{float}
 \usepackage{multicol}
 \usepackage[all]{xy}

\newtheorem{proposition}{Proposition}[section]
\newtheorem{theoremb}[proposition]{Theorem}

\theoremstyle{definition}

\theoremstyle{remark}

\newcommand{\be}{\begin{equation}}
\newcommand{\ee}{\end{equation}}

\newcommand{\bd}{\begin{displaymath}}
\newcommand{\ed}{\end{displaymath}}
\newcommand{\ben}{\begin{enumerate}}
\newcommand{\een}{\end{enumerate}}

\title{Rational homotopy theory via Sullivan models: A survey} 
\author{Yves F\'elix and Steve Halperin}
\date{\today}

\begin{document}

\maketitle

 Rational homotopy theory begins with Sullivan's introduction of the localization of a topological space (\cite{Su1})  although as he states in his preface, the \emph{compulsion to localize began with the author's work on invariants of combinatorial manifolds in 1965-67.} The two foundational articles (\cite{S}, \cite{Q}) then firmly established the subject as an independent sub-discipline of algebraic topology. 
 
The two distinct approaches by Quillen and Sullivan directly reflect the original two means of computing homology and cohomology at the start of the twentieth century: simplicial homology of a polyhedron, and deRham cohomology of a manifold.  On the one hand, Quillen establishes an equivalence of homotopy categories between simply connected rational spaces and a class of differential graded Lie algebras (dgl's). On the other, inspired by the differential forms on a manifold, Sullivan constructs a morphism $A_{PL}$ of homotopy categories from all topological spaces, $X$, to rational commutative differential cochain algebras (cdga's).

 For simply connected rational spaces with finite Betti numbers, and the corresponding category of cdga's, this is also an equivalence. (Hess \cite{Hhist} provides a more detailed and excellent history of the subject.) 
Moreover, if a cdga $(A,d)$ corresponds to a path connected space $X$, then with it Sullivan associates its minimal model, a cdga which is free as a commutative graded algebra and whose isomorphism class depends only on the homotopy type of $X$. These minimal models, which can often be computed, provide detailed information about $X$ and its invariants, including its Lusternik-Schirelmann category and a new invariant depth$\, X$. In particular, properties of fundamental groups, including their Malcev completions, are accessible with this approach. 

Our objective, then, is to provide an overview of rational homotopy theory as viewed through the lens of Sullivan's minimal models.  We describe many of the main results, and a range of examples and applications new and old. However, for ease of reading the proofs  and what are now classical constructions  are not included. For these the reader is generically referred to \cite{FHT} and \cite{FHTII} or to the excellent introductions to the subject by Hess \cite{Hsu} and Berglund \cite{Be1}; other references are provided in the text where required. We also include a range of open problems, new and old, distributed throughout the text.

The breadth of established material, as well as ongoing new results, applications and conjectures, poses a challenge to authors of such a survey, which cannot hope to be encyclopaedic. We have had to be selective, understanding that others might well have made different choices.  

Finally, we would like to express our profound appreciation to the other rational homotopy theorists with whom we have had the good fortune to be able to collaborate over the last 40 years, and most especially to our good friend, Jean-Claude Thomas, frequently our co-author, including in particular for the two volumes on rational homotopy.

\tableofcontents

\section{Notation and terminology}
 
Unless otherwise specified we work over $\mathbb Q$ as ground field, noting that for much of the material any field of characteristic zero could be used.

A \emph{graded vector space} $S$ is a collection $\{S_p\}_{p\in \mathbb Z}$ or a collection $\{S^p\}_{p\in \mathbb Z}$ of vector spaces, with the convention $S_p=S^{-p}$ used to avoid negative indices. A graded vector space $S$ has \emph{finite type} if dim$\, S^p<\infty$, $p\in \mathbb Z$.

If $S$ and $T$ are graded vector spaces then Hom$(S,T)$ is the graded space of \emph{linear maps} given by Hom$_k(S,T) = \prod_{r\in \mathbb Z} \mbox{Hom}(S_r, T_{r+k})$, but \emph{morphisms} are required to preserve degrees. For simplicity, Hom$(S, \mathbb Q)$ is then denoted by $S^\#$. When $S= \{S^p\}$ an  element $x\in S^p$ has \emph{degree} $p$.

The \emph{suspension} $sS$ is the graded vector space defined by $(sS)_p = S_{p-1}$; this isomorphism is denoted $x\mapsto sx$. The inverse isomorphism is denoted by $s^{-1}$. 

\emph{Graded algebras} $A$ are graded vector spaces together with an associative multiplication $A\otimes A\to A$ of degree zero, and have an identity $1\in A_0$. A graded algebra is \emph{commutative} (cga for short) if for $a,b\in A$,
$$ab = (-1)^{deg\, a \cdot deg\, b} ba\,.$$

\emph{Graded Lie algebras}, $L$, are graded vector spaces equipped with a Lie bracket $[\hspace{2mm}, \hspace{2mm}] : L\otimes L\to L$ of degree $0$ and satisfying for $x,y,z\in L$:
$$[x,y] + (-1)^{deg\, x\cdot deg\, y}[y,x]= 0\hspace{5mm} \mbox{and } [x,[y,z]] = [[x,y],z] + (-1)^{deg\, x\cdot deg\, y} [y, [x,z]]\,.$$

The \emph{universal enveloping algebra}, $UL$, is the graded algebra $TL/J$, where $TL$ is the tensor algebra and $J$ is the ideal generated by the elements $[x,y] - (x\otimes y - (-1)^{deg\, x\cdot deg\, y} y\otimes x)$.

A \emph{complex}  or a \emph{differential graded vector space}, $(S,d)$, is a graded vector space together with a linear map $d\in \mbox{Hom}_{-1}(S,S)$ satisfying $d^2= 0$. Its \emph{homology} $H(S)$ is the graded vector space ker$\, d /$ Im$\, d$. A \emph{quasi-isomorphism} in a category of complexes, denoted by $\stackrel{\simeq}{\longrightarrow}$, is a morphism $\varphi$ satisfying $H(\varphi)$ is an isomorphism. A \emph{differential graded algebra} (dga for short) is a pair $(A,d)$ in which $A$ is a graded algebra, $(A,d)$ is a complex,  and $d$ is a \emph{derivation}: $d(ab)= (da)b + (-1){deg\, a} a (db)$. If $A$ is commutative then $A$ is a cdga. A \emph{differential graded Lie algebra} (dgl for short) is a pair $(L,\partial)$ in which $L$ is a graded Lie algebra, $(L, \partial)$ is a complex, and $\partial$ is a derivation: $\partial [x,y]= [\partial x, y] + (-1)^{deg\, x} [x, \partial y]$. 

If $(A,d)$ is a dga then an $(A,d)$-module is a pair $(M,d)$ in which $M$ is a (graded) $A$-module, $(M,d)$ is a complex, and $d(a\cdot m) = (da)\cdot m + (-1)^{deg\, a} a\cdot dm$. An $(A,d)$-module, $(P,d)$, is \emph{semi-free} if $P$ is the union of an increasing sequence $(P_0,d)\subset \cdots \subset (P_r,d)\subset \cdots $  of sub $(A,d)$-modules in which $(P_0,d)$ and each $(P_{r+1}/P_r, d)$ are direct sums of copies of $(A,d)$. Any $(A,d)$-module $M,d)$ admits a quasi-isomorphism
$$\varphi : (P,d)\stackrel{\simeq}{\longrightarrow} (M,d)$$
from a semi-free $(A,d)$-module.

Finally, for any topological space $X$, the rational cohomology $H^*(X;\mathbb Q)$ is denoted simply by $H(X)$, to simplify the notation in the many places this appears. Otherwise homology is always denoted with the coefficient ring.

 \section{Sullivan models}
 
 Sullivan's approach to rational homotopy theory is accomplished by a passage from topological spaces $X$ to commutative differential graded algebras (cdga's) via a quasi-isomorphism
 $$\varphi : (\land V,d) \to A_{PL}(X)\,,$$
 where in particular
 \begin{enumerate}\item[(i)] $\land V$ is the free commutative graded algebra on a graded vector space $V = \{V^n\}_{n\geq 1}$; $\land V = \oplus_{r\geq 0} \land^rV$ where $\land^rV$ is the linear span of the monomials in $V$ of length $r$.
 \item[(ii)] $\xymatrix{A_{PL}: \mbox{Top}\ar@{~>}[r] & \mbox{Cdga}}$   is a contravariant functor for which the graded algebras $H(X)$ and $H(A_{PL}(X))$ are naturally isomorphic.
 \end{enumerate}
  This construction has three steps. 
  
\vspace{2mm}{\bf The first step} is the definition of the simplicial cdga $(A_{PL})_*$ and the definition of $A_{PL}$ as the composite of two functors
$$\xymatrix{\mbox{Top} \ar@{~>}[r]^{Sing} & \mbox{Simpl} \ar@{~>}[r]^{A_{PL}} & \mbox{Cdga}\,,}$$
where Simpl is the category of simplicial sets, $Sing (X)$ is the simplicial set of singular simplices in $X$, and $A_{PL}(S) = \mbox{Simpl}(S, (A_{PL})_*)$. Here $(A_{PL})_*$ is defined as follows: let $t_i$ denote variables of degree $0$, then
$$(A_{PL})_k = \land (t_0, \cdots t_k, dt_0, \cdots , dt_k)/ (\sum t_i)  - 1, \sum dt_i\,.$$
The face and degeneracy morphisms are given by
$$\partial_it_j = \left\{\begin{array}{ll} t_j\,, & j<i\\0\,, & j= i\\t_{j-1}\,, &j>i\end{array}\right. \hspace{1cm}\mbox{and } s_jt_i = \left\{\begin{array}{ll} t_i\,, & i<j\\
t_{i+1}\,, & i\geq j\,.\end{array}\right.$$
Thus $\varphi\in A_{PL}^n(X):=A_{PL}( Sing \, X)$ assigns to each singular $k$ simplex in $X$ an element of degree $n$ in $(A_{PL})_k$, compatibly with the face and degeneracy maps.

\vspace{2mm}{\bf The second step} is the introduction of a distinguished class of cdga's, the Sullivan algebras, and the establishment of their properties.

\vspace{3mm}\noindent {\bf Definition.} \begin{enumerate}
\item[(i)] A \emph{Sullivan algebra} is a cdga of the form $(\land V,d)$ in which $V = \{V^k\}_{k\geq 1}$ is the union of subspaces $V(0)\subset \cdots \subset V(r)\subset \cdots $ in which $d: V(0)\to 0$ and $d : V(r+1)\to \land V(r)$ for $r\geq 0$. (Note : a Sullivan algebra will admit many choices of generating space $V$ !).
\item[(ii)] $(\land V,d)$ is \emph{minimal}  if $d: V\to \land^{\geq 2}V$, and \emph{contractible} if for some $W\subset V$ the inclusions of $W$ and $d(W)$ in $\land V$ extend to an isomorphism
$$\land W \otimes \land d(W)\stackrel{\cong}{\to} \land V\,.$$ 
\end{enumerate}

Arguments from linear algebra then establish the

\vspace{3mm}\noindent {\bf Properties}
\begin{enumerate}
\item[(i)] A Sullivan algebra is the tensor product of a minimal algebra and a contractible Sullivan algebra

\item[(ii)] A Sullivan algebra $(\land V,d)$ is the direct limit of the sub Sullivan algebras $(\land V_\alpha,d)$ in which $V_\alpha\subset V$ and dim$\, V_\alpha<\infty$.
\item[(iii)] A quasi-isomorphism between minimal Sullivan algebras is an isomorphism.
\end{enumerate}

Additionally, in analogy with the inclusion of the endpoints in an interval, the cdga $\land (t,dt)$ with deg$\, t= 0$ has augmentations $\varepsilon_0, \varepsilon_1 : \land (t,dt)\to \mathbb Q$ defined by $\varepsilon_0(t) = 0$ and $\varepsilon_1(t) = 1$.

\vspace{3mm}\noindent {\bf Definition.} Two cdga morphisms $\varphi_0, \varphi_1 : (\land V,d)\to (A,d)$ from a Sullivan algebra are \emph{homotopic} via a homotopy $\Phi : (\land V,d)\to \land (t,dt)\otimes (A,d)$ if $(\varepsilon_i\otimes id) \circ \Phi = \varphi_i$. This is denoted $\varphi_0\sim \varphi_1$. 

\vspace{3mm}Moreover, the cdga analogue of a based topological space $X$ is an augmented cdga $\varepsilon : (A,d)\to \mathbb Q$: if $j :x\to X$ is the inclusion of a base point then $A_{PL}(j) : A_{PL}(X)\to A_{PL}(x)= \mathbb Q$ is an augmentation. Thus if in the definition above $\varepsilon : (A,d)\to \mathbb Q$ is an augmentation, then $\Phi$ is a \emph{based homotopy} if $\Phi(\land^+V)\subset \land (t_0, t_1) \otimes \mbox{ker}\, \varepsilon$, and this is denoted $\varphi_0\sim_*\varphi_1$. 

As in the topological category, homotopy and based homotopy are equivalence relations.

Finally, the essential property of Sullivan algebras $(\land V,d)$ is   that morphisms from Sullivan algebras lift up to homotopy through quasi-isomorphisms. More precisely, given $\varphi$ and $\eta$ in the diagram below
$$\xymatrix{\mbox{} & (C,d)\ar[d]_\simeq^\eta\\
(\land V,d) \ar@{-->}[ru]^\psi \ar[r]_\varphi & (A,d)\,,}$$
there is a unique homotopy class of morphisms $\psi$ such that $\eta\circ \psi \sim \varphi$. Moreover, the homotopy class only depends on the homotopy class of $\varphi$. The analogous result also holds for based homotopy when $\eta$ is a morphism of augmented cdga's.

\vspace{2mm}{\bf The third and final step} is the

\begin{theoremb}\label{minmod} For each path connected topological space $X$ there is a quasi-isomorphism
$$\varphi : (\land V,d)\stackrel{\simeq}{\to} A_{PL}(X)$$
from a minimal Sullivan algebra. If $\psi : (\land W,d)\stackrel{\simeq}{\to} A_{PL}(X)$ is a second quasi-isomorphism, then there is an isomorphism
$$\chi : (\land W,d)\stackrel{\cong}{\to} (\land V,d)$$
such that $\varphi\circ \chi \sim \psi$. If $X$ is a based space then $\chi$ may be chosen so that $\varphi\circ \chi\sim_* \varphi$.\end{theoremb}

 \noindent {\bf Definition.} The morphism $\varphi : (\land V,d)\stackrel{\simeq}{\to} A_{PL}(X)$ is a \emph{minimal Sullivan model} for $X$.

\vspace{3mm} Similarly, given a map $f : X\to Y$ of (based) topological spaces, minimal Sullivan models for $X$ and $Y$ embed in a (based) homotopy commutative diagram
$$
\xymatrix{ (\land V,d) \ar[r]^\simeq & A_{PL}(X)\\
(\land W,d) \ar[u]^\varphi \ar[r]^\simeq & A_{PL}(Y) \ar[u]_{A_{PL}(f)}
}$$
in which the (based) homotopy class of $\varphi$ is determined by the (based) homotopy class of $f$. The morphism $\varphi$ is a \emph{Sullivan representative} of $f$.

\vspace{3mm}\noindent {\bf Remarks 1.} Theorem 1 remains true if $A_{PL}(X)$ is replaced by any cdga $(A,d)$ for which $H^0(A,d) = \mathbb Q$. In this case a quasi-isomorphism $(\land V,d) \stackrel{\simeq}{\to} (A,d)$ from a minimal Sullivan algebra is a \emph{minimal Sullivan model} for $(A,d)$.

\vspace{2mm}{\bf 2.} Two cdga's $(A,d)$ and $(B,d)$ are \emph{weakly equivalent} if there is a finite chain of cdga quasi-isomorphisms
$$(A,d) = (C(0),d)\stackrel{\simeq}{\to} (C(1),d) \stackrel{\simeq}{\leftarrow} \cdots \stackrel{\simeq}{\to} (C(p),d) = (B,d)\,.$$
In this case we write $(A,d)\simeq (B,d)$. If $H^0(A,d) =\mathbb Q = H^0(B,d)$ then $(A,d)\simeq (B,d)$ if and only if for some minimal Sullivan algebra $(\land V,d)$ there are quasi-isomorphisms
$$(A,d)\stackrel{\simeq}{\leftarrow} (\land V,d)\stackrel{\simeq}{\to} (B,d)\,.$$

The cdga $A_{PL}(X)$, and by consequence the minimal Sullivan model of $X$ is directly related to other dga's associated with $X$:
\begin{enumerate}
\item[1.] Denote by $C^*(X;\mathbb Q)$ the usual cochain algebra on $X$ with rational coefficients, then there is a sequence of quasi-isomorphisms of dga's
$$A_{PL}(X) \stackrel{\simeq}{\longleftarrow} E(X) \stackrel{\simeq}{\longrightarrow} C^*(X)\,.$$
Two dga's $(A,d)$ and $(B,d)$ are \emph{equivalent} if they are connected by quasi-isomorphisms,
$$\xymatrix{(A,d)\ar[r]^\simeq & (A(1),d) & \cdots \ar[l]_\simeq \ar[r]^\simeq & (A(k),d) & (B,d)\ar[l]_\simeq.}$$
When a simply connected CW complex $X$ has finite Betti numbers then Adams-Hilton (\cite{AH}) construct a natural bijection between the equivalence classes of $C^*(X;\mathbb Q)$ and $C_*(\Omega X;\mathbb Q)$, $\Omega X$ denoting the loop space. 

\hspace{3mm}
It is an open question to construct two spaces $X$ and $Y$ with equivalent cochain algebras, $C^*(X;\mathbb Q) \simeq C^*(Y;\mathbb Q)$, but different minimal Sullivan models. Is there an obstruction theory to measure the difference ?
\item[2.] If $X$ is a   connected manifold, then $A_{PL}(X)\otimes_{\mathbb Q}\mathbb R$ is weakly equivalent to the cdga of de Rham forms on $X$, $A_{PL}(X)\otimes_{\mathbb Q}\mathbb R \simeq A_{DR}(X)$. 
\item[3.] It is an open question whether geometric data (as opposed to triangulations) on a manifold $M$ could permit the construction of a rational Sullivan model. However (\cite{Hal}) a Riemannian metric does determine an algorithm for computing an finitely generated field extension, $l\!k$ of $\mathbb Q$ and the minimal Sullivan model defined over $l\!k$. 
\item[4.] The existence of a free torus action on a finite CW complex $X$ imposes restrictions on the minimal Sullivan model. Could these be extended to resolve the question: in this case is dim$\, H(X)\geq 2^n$ ?
\end{enumerate}

\vspace{3mm}\noindent {\bf Examples of Sullivan models}
\begin{enumerate}
\item Denote by $M_X$ the minimal Sullivan model of a space $X$. Then
$$\renewcommand{\arraystretch}{1.3}\begin{array}{l}
M_{S^{2n+1} } = (\land u, 0)\,, \hspace{1cm} \mbox{deg}\, u = 2n+1\\
 M_{S^{2n}} = (\land (a,b),d)\,, \hspace{5mm} da=0, db=a^2\,, \hspace{1cm} \mbox{deg}\, a= 2n\\
M_{X\times Y} \cong M_X\otimes M_Y\,, \mbox{if one of $H(X)$ or $H(Y)$ is a graded vector space of finite type}\\
M_{X\vee Y} \simeq M_X\oplus_{\mathbb Q} M_Y\\
M_{K(\mathbb Z,n)} = (\land a, 0)\,, \hspace{1cm} \mbox{deg}\, a= n
\end{array}
\renewcommand{\arraystretch}{1}
$$
\item \label{hospace} The minimal Sullivan model of a connected compact Lie group is a cdga of the form $(\land V_G,0)$ where $V_G$ is finite dimensional and concentrated in odd degrees. A Sullivan model of its classifying space $BG$ is $(\land s^{-1}V_G,0)$ with as usual $(s^{-1}V_G)^n = V_G^{n-1}$. If $X= G/H$ is a  homogeneous space where $G$ and $H$ are compact connected Lie groups, denote by $\varphi : (\land s^{-1}V_G,0)\to (\land s^{-1}V_H,0)$ a Sullivan representative of the induced map $f : BH\to BG$. Then a minimal Sullivan model for $G/H$ is 
$$(\land s^{-1}V_H\otimes \land V_G,d)$$
where $d(s^{-1}V_H)= 0$ and $d(x) = \varphi (s^{-1}x)$ for $x\in V_G$.
\item Suppose a compact connected Lie group $G$ acts on a manifold $M$, then the injection of the algebra of $G$-invariants forms $A_{DR}^GM \to A_{DR}(M)$ is a quasi-isomorphism.
\item Let $f : X\to Y$ be a continuous map between simply connected spaces and let $\varphi : (\land V,d)\to (\land V,d)$ be a representative of $f$. Extend $\varphi$ to a surjective map $\overline{\varphi} : (\land W,d)\otimes \land (S\oplus dS)\to (\land V,d)$, from the tensor product with a contractible Sullivan algebra. Then a minimal Sullivan model for $\mathbb Q\oplus \mbox{ker}\, \overline{\varphi}$ is a minimal Sullivan model for $C_f$, the homotopy cofibre of $f$.
\item A space $X$ with minimal Sullivan model $(\land V,d)$ is \emph{formal} if there is a quasi-isomorphism $\varphi : (\land V,d)\to (H(X),0)$. Every $(n-1)$-connected space, $n\geq  2$, of dimension $\leq 3n-2$ is formal. The $H$-spaces,   symmetric spaces of compact type, and   compact K\"ahler manifolds (\cite{DGMS}) are formal. In a formal space all Massey products of order $\geq 3$ are trivial, and so   a non-trivial Massey product is   an obstruction to formality. In \cite{HaS} Halperin and Stasheff have given a simple algorithm for deciding the formality of a space. 

There are many non formal spaces. For instance $(\land u,v,w),d)$ with $du=dv=0$ and $dw= uv$ with $u$ and $v$ in odd degrees is a non formal space. In geometry \cite[Sec.10.12]{GHV} gives a Lie theoretic characterization of formal homogeneous spaces, while  examples of non-formal simply connected symplectic manifolds have been given by Babenko and Taimanov (\cite{BT}). Moreover, in \cite{CFM} Cavalcanti, Fern\'andez and Munoz give an example of a non formal simply-connected compact symplectic manifold of dimension 8 which satisfies the Lefchetz property,  raising the question whether there is  a nice generalization of  K\"ahler that implies formality. 
\end{enumerate}

\section{The spatial realization of a Sullivan algebra}

The second component of Sullivan's approach to rational homotopy theory is the passage from Sullivan algebras to CW complexes via the composite of two functors
$$ \xymatrix{\vert \,\,\vert : \mbox{Cdga}\ar@{~>}[r] &\mbox{Simpl} \ar@{~>}[r] & \mbox{CW}}$$
respectively adjoint to $A_{PL}$ and to $Sing$. (As with $A_{PL}$, we use the same notation for the composite and for the second functor alone.)

Here
 $<A,d>_n = \mbox{Cdga} ((A,d), (A_{PL})_n)\,,$ 
while $\vert\,\,\vert $ is the Milnor realization.

In particular, if $X$ is a CW complex then adjoint to $id_{Sing\, X}$ is a (based) homotopy equivalence $\vert Sing\, X\vert \to X$ which determines up to (based) homotopy a (based) homotopy equivalence
$$r_X : \vert Sing\, X\vert \to X\,.$$

\vspace{3mm}\noindent {\bf Definition.} The \emph{spatial realization} of a Sullivan algebra $(\land V,d) $ is the based CW complex $\vert \land V,d\vert: = \vert <\land V,d>\vert$. (Note that $\vert \land V,d\vert$ has a single $0$-cell which serves as the base point.)

\vspace{3mm}There are two distinct contexts in which it is useful to convert a cdga morphism to a continuous map. Both constructions preserve (based) homotopy.

\vspace{2mm}
(i) If $\varphi : (\land V,d) \to (\land W,d)$ is a morphism of Sullivan algebras, then $\varphi\mapsto \vert \varphi\vert$,
$$\vert \varphi \vert : \vert \land V,d\vert \leftarrow \vert \land W,d\vert,$$
is   the functor described above, applied to morphisms.

\vspace{2mm}
(ii) If $\varphi : (\land V,d) \to A_{PL}(X)$ is a morphism from a Sullivan algebra and if $X$ is a CW complex, then by adjunction, $\varphi$ induces a map
$$\widehat{\varphi} : X\to \vert \land V,d\vert\,.$$

\vspace{3mm}\noindent {\bf Example.}  Suppose $\varphi : (\land V,d)\to A_{PL}(X)$ is the minimal Sullivan model of a CW complex with finite Betti numbers. If $H^1(X)= 0$ then $V = V^{\geq 2}$ is a graded vector space of finite type, and 
$$H(\widehat{\varphi}) : H(\vert\land V,d\vert)\stackrel{\cong}{\longrightarrow} H(X)$$
is an isomorphism.
If, in addition, $X$ is simply connected then 
$$\pi_*(\widehat{\varphi}) : \pi_*(X)\otimes \mathbb Q \to \pi_*(\vert \land V,d\vert)$$
is also an isomorphism. In this case $\widehat{\varphi} : X\to \vert \land V,d\vert$ is the \emph{rationalization} of $X$.

\vspace{5mm} This example generalizes under certain conditions to non-simply connected spaces. For instance Sullivan's infinite telescope is a rationalization $S^1\to S^1_{\mathbb Q}$ of the circle. But, as the next example shows, when $X$ is not simply connected $\pi_*\vert \land V,d\vert$ may be quite different from $\pi_*(X)\otimes \mathbb Q$ in degrees $\geq 2$ even when $H^1(X)= 0$ and the Betti numbers of $X$ are finite. In fact the determination of the cohomology of $\vert\land V,d\vert$ when $V^1$ is infinite dimensional remains an open question.

\vspace{3mm}\noindent {\bf Example.} The involution $\tau : (x,y,z)\to (-x,-y,-z)$ in $S^2\times S^2\times S^2$ acts freely, and so dividing by $\tau$ gives a CW complex $X$. The cohomology algebra of $X$ is the sub algebra of $H(S^2\times S^2\times S^2)$ of classes left fixed by $\tau$. Let $\alpha, \beta, \gamma$ be the fundamental classes for the three $2$-spheres. Then $1, \alpha\beta, \alpha\gamma $ and $\beta\gamma$ are a basis of $H(X)$. It follows that $H(X) = H(S^4\vee S^4\vee S^4)$.

But $S^4\vee S^4\vee S^4$ is intrinsically formal, which implies that the minimal Sullivan model $(\land V,d)$ of $X$ is also the minimal Sullivan model of $S^4\vee S^4\vee S^4$ and so
$$\vert \land V,d\vert = (S^4\vee S^4\vee S^4)_{\mathbb Q}\,.$$
A simple computation also shows that in this case, while each dim$\, V^k<\infty$, $V$ itself is infinite dimensional.
Moreover, while   $\pi_{\geq 2}(X)\otimes \mathbb Q$ is concentrated in degrees $2$ and $3$, $\pi_*\vert \land V\vert$ is infinite dimensional and concentrated in degrees $3n+1$.

\vspace{3mm}
Finally, for each Sullivan algebra $(\land V,d)$, adjoint to $id_{(\land V,d)}$ is a morphism, natural in $(\land V,d)$, 
$$m_V : (\land V,d) \to A_{PL}\vert \land V,d\vert\,,\hspace{1cm}\mbox{with }  \vert m_V\vert = id_{\vert \land V,d\vert}\,.$$
 In general it may not be a quasi-isomorphism, but if $V$ is a graded vector space of finite type then $H(m_V)$ is an isomorphism and $m_V$ is the minimal Sullivan model.  

Further, for any connected CW complex $X$ and any morphism $\varphi : (\land V,d)\to A_{PL}(X)$ from a Sullivan algebra,
$$\varphi \sim_* A_{PL}(\widehat{\varphi}) \circ m_V\,.$$
In particular, if $\varphi$ is the minimal Sullivan model of $X$ then $H(m_V)$ is injective and $H\,(\widehat{\varphi})$ is surjective. Thus  $H(m_V)$ is an isomorphism if and only if $H\,(\widehat{\varphi})$ is an isomorphism.

\section{Homotopy groups}

Associated with any minimal Sullivan algebra $(\land V,d)$ are its homotopy groups $\pi_*(\land V,d)$, whose definition is motivated by the following construction.

Suppose $X$ is a based connected CW complex and 
$$\varphi : (\land V,d)\to A_{PL}(X)$$ is a morphism. Then $\varphi$ determines a set map
$$\pi_*(\varphi) : \pi_*(X)\to \left( \land^+V/\land^+V{\scriptstyle \land} \land^+V\right)^\#= V^\#$$
as follows:
\begin{enumerate}
\item[$\bullet$] Represent $\alpha\in \pi_n(X)$ by $f : S^n\to X$, so that $A_{PL}(f)\circ \varphi : (\land V,d)\to A_{PL}(S^n)$. Since $A_{PL}(S^n)\simeq H(S^n)$, this morphism lifts (up to based homotopy) to a morphism
$$\varphi_f : (\land V,d) \to H(S^n)\,.$$
\item[$\bullet$] The standard orientation in $\Delta^n$ together with the map $\Delta^n \to \Delta^n/\partial \Delta^n= S^n$ defines an orientation class $\omega_n \in H^n(S^n)$.
\item[$\bullet$] Define $\pi_*(\varphi)\alpha$ by
$$(\pi_*(\varphi)\alpha)(\Phi)\cdot \omega_n = H^n(\varphi_f)(\Phi)\,, \hspace{1cm} \Phi\in V^n\,.$$
\end{enumerate}

\vspace{3mm}\noindent {\bf Remark.} The set map $\pi_*(\varphi)$ is natural with respect to morphisms of minimal Sullivan algebras and with respect to based maps of CW complexes.

\begin{theoremb} When $\varphi : (\land V,d)\to A_{PL}(X)$ is the minimal Sullivan model of a simply connected space with finite Betti numbers, then $\pi_*(\varphi) : \pi_*(X)\otimes \mathbb Q\to V^\#$ is a linear isomorphism.\end{theoremb}

\vspace{3mm} The construction above applies in particular to the morphism $$m_V : (\land V,d) \to A_{PL}\vert \land V,d\vert\,,$$ where we have the

\begin{theoremb} For any minimal Sullivan algebra $(\land V,d)$,
\begin{enumerate}
\item[(i)] $\pi_k(m_V): \pi_k\vert \land V,d\vert \to (V^k)^{\#}$ is a linear isomorphism for $k\geq 2$.
\item[(ii)] $\pi_1(m_V): \pi_1\vert \land V,d\vert \to (V^1)^\#$ is a bijection.
\end{enumerate}\end{theoremb}

\noindent \emph{In particular, when $k= 1$,  $\pi_1(m_V)$ induces a group structure on $(V^1)^\#$}:

\vspace{3mm}\noindent {\bf Definition.} The groups $\pi_k(\land V,d): = (V^k)^\#$, $k\geq 1$ are the \emph{homotopy groups} of the minimal Sullivan algebra $(\land V,d)$.

\vspace{3mm}\noindent {\bf Remarks 1.}  The multiplication in $\pi_1(\land V,d)$ will be made explicit in the next section.

\vspace{2mm}{\bf 2.} By naturality it follows that if $X$ is a CW complex and $\varphi : (\land V,d)\to A_{PL}(X)$ is a morphism from a minimal Sullivan algebra then the commutative diagram
$$\xymatrix{
\mbox{} && \pi_*\vert \land V,d\vert \ar[dd]^\cong_{\pi_* (m_V)}\\
\pi_*(X) \ar[rru]^{\pi_*\widehat{\varphi}}
\ar[rrd]_{\pi_*(\varphi)} \\
\mbox{} && V^\#}$$
identifies $\pi_*(\varphi)$ with the homomorphism $\pi_*\widehat{\varphi}$.

\vspace{2mm} {\bf 3.} When $\varphi : (\land V,d) \to A_{PL}(X)$ is the minimal Sullivan model of a based connected CW complex then
$$H(\land V,d) \cong H(X) \hspace{1cm}\mbox{and } \pi_*(\land V,d) \cong \pi_*\vert \land V,d\vert\,.$$
However it may happen that neither $H(\widehat{\varphi})$ nor $\pi_*\widehat{\varphi}$ is an isomorphism.

\vspace{3mm}\noindent {\bf Example.} Let $(\land Z,d)$ be the minimal model of the sphere $S^n$, for some $n\geq 1$. Since $S^n$ is formal, there is a quasi-isomorphism $\rho : (\land Z,d)\to (H(S^n),0)\cong (\land u/u^2, 0)$.   Next, let $g : S^n\times X\to Y$ be a continuous map and $f$ its restriction to $\{x_0\}\times X$ where $x_0$ is a base point in $S^n$. The Sullivan representative for $g$  (followed by $\rho\otimes id$) has the form 
$$\psi : (\land V,d)\to (\land u/u^2, 0\otimes \land W,d)\,, \hspace{1cm} \psi (v)= \varphi (v) + u\theta (v)\,, \hspace{5mm} v\in V\,.$$
Here $\varphi$ is a Sullivan representative for $f$, and $\theta : \land V\to \land W$ is a $\varphi$-derivation, i.e. a linear map satisfying
$$\theta (ab)= \theta (a)\varphi (b) + (-1)^{deg\, a\cdot deg\, \theta} \varphi (a)\theta (b)\,.$$

Suppose $X$ is a finite CW complex and $Y$ is a nilpotent space with finite Betti numbers. Denote by $F(X,Y)$ the mapping space of continuous maps from $X$ to $Y$, and by $F(X,Y;f)$ the path component of $g$. Then    $F(X,Y;f)$ is a nilpotent space (\cite{Hil}), and the above correspondence $g\mapsto \theta$ induces for $n\geq 2$ an isomorphism, and for $n=1$ a bijection,
$$\pi_n(F(X,Y;f))\otimes \mathbb Q \stackrel{\cong}{\longrightarrow} H_{n}\left(\mbox{Der}_\varphi (\land V, \land W), D)\right)\,,$$
where Der$_\varphi$ denotes the vector space of $\varphi$-derivations and $D\theta = d\theta - (-1)^{deg\, \theta} \theta d$ (\cite{BM}, \cite{LuS}).

in particular, if $f$ is the constant map, then $\mbox{Der}_\varphi (\land V, \land W)= \mbox{Hom}(V, \land W)$ and 
$$\pi_n(F(X,Y;f))\otimes \mathbb Q \cong \oplus_q \mbox{Hom}(\pi_q(Y), H_{q-n}(X))\,.$$

 \section{The homotopy Lie algebra and the fundamental group}
 
Fix a minimal Sullivan algebra $(\land V,d)$.  Theorem 4.2 identifies the abelian groups $\pi_k\vert\land V,d\vert$, $k\geq 2$, as the rational vector space $(V^k)^\#$. To make explicit the multiplication in $\pi_1(\land V,d)$ we introduce the homotopy Lie algebra of $(\land V,d)$, defined next.  

First we establish the convention that in the pairing between $V$ and $V^\#$,  $V^\#$ will act from the right so that the pairing is written
$$V\times V^\#\to \mathbb Q\,, \hspace{1cm} v,f\mapsto <v,f>\,.$$
This induces the pairing $\land^2V \times (V^\#\times V^\#)\to \mathbb Q$ given by
$$<v\land w,f,g> = <v,g>\, <w, f> + (-1)^{deg\, v\cdot deg\, w} <v,f>\, <w,g>.$$

Now denote by $d_1v$ the component of $dv$ in $\land^2V$, $(\land V, d_1)$ is a \emph{quadratic Sullivan algebra}.
 
 \vspace{3mm}\noindent {\bf Definition.} The \emph{homotopy Lie algebra}, $L_V$ is the graded Lie algebra in which $sL_V$ is the graded vector space $V^\#$, and 
 $$(v, s[x,y]> = (-1)^{deg\, y+1}<d_1v, sx, sy>\,.$$
 
\vspace{1mm}\noindent Thus If $V^1\neq 0$ then $(L_V)_0\neq 0$.

\vspace{2mm} For simplicity, in the rest of this section, $L_V$ will be denoted simply by $L$.

 \vspace{3mm}
To describe $\pi_1(\land V,d)$, let $I^n$ denote the $n^{th}$ power of the augmentation ideal in the universal enveloping algebra $UL$. Then $\widehat{UL}:= \varprojlim_n UL/I^n$ is the \emph{completion} of $UL$. In particular, an injective set map $L_0\stackrel{\exp}{\longrightarrow} \widehat{UL}$ is given by 
   $\exp (x) = \sum_{n=0}^\infty x^n/n!\,.$ 
  
 \begin{theoremb} Suppose dim$\, H^1(\land V,d)<\infty$. Then multiplication in $\widehat{UL}$ restricts to a product in $\exp (L_0)$ which makes $G_L:= \exp (L_0)$ into a group, and
 $$\exp\circ s^{-1} : \pi_1(\land V,d) \stackrel{\cong}{\rightarrow} G_L$$
 is an isomorphism of groups. \end{theoremb}
 
 \vspace{3mm}\noindent {\bf Remarks 1.} The 
  inverse of $\exp$ is the formal power series   $\log: G_L\to L_0$,  
 $\log a = \sum_{n\geq 1} (-1)^{n-1} \,\,\frac{(a-1)^n}{n}\,.$ Thus the product in $\pi_1(\land V,d)$ is given explicitly by
 $$\alpha\cdot \beta = s \log (\exp s^{-1}\alpha\cdot \exp s^{-1}\beta)\,.$$

 \vspace{2mm}{\bf 2.}  While the theorem has the restriction dim$\, H^1(\land V,d)<\infty$,  it does apply when dim$\,V<\infty$. Since $(\land V,d) = \varinjlim_\alpha (\land V_\alpha,d)$ with dim$\, V_\alpha<\infty$, it follows that in general
 $$\pi_1(\land V,d)= \varprojlim_\alpha G_\alpha$$
 where $G_\alpha$ if the "exponential group" of $(\land V_\alpha,d)$.

 \vspace{3mm}
Associated with the group $\pi_1 :=\pi_1(\land V,d)$ and the Lie algebra $L_0$ are their respective  lower central series $\pi_1= \pi_1^1\supset \cdots \supset \pi_1^r\supset \cdots $ and $L_0= L_0^1\supset \cdots \supset L_0^r\supset \cdots $. Here $\pi_1^{r+1}$ and $L_0^{r+1}$ are respectively the subgroup generated by the commutators $[\alpha, \beta]$, $\alpha\in \pi_1$, $\beta\in \pi_1^{r}$ and the linear span of the commutators $[x,y]$, $x\in L_0$, $y\in L_0^r$. 

On the other hand, $d_1$ determines the increasing filtration $V_0^1\subset \cdots V_r^1\subset \cdots$ given by 
$$V_{0}^1 = V^1\cap \mbox{ker}\, d_1\hspace{5mm}\mbox{ and } V_{r+1}^1= d_1^{-1} (\land V_{r}^1)\,.$$
The defining condition for Sullivan algebras gives $V^1 = \cup_m V_{m}^1$. In particular, $$V_0^1 = H^1(\land V,d)\,.$$   These three filtrations are closely related when dim$\, H^1(\land V,d)<\infty$: 

 \begin{theoremb} The vector spaces $\pi_1/\pi_1^2, L_0/L_0^2,$ and $(V_0^1)^\#$ are isomorphic if one of them is finite dimensional. In this case
 \begin{enumerate}
 \item[(i)] $s: L_0\stackrel{\cong}{\rightarrow } \pi_1$ restricts to bijections $s: L_0^m\stackrel{\cong}{\rightarrow} \pi_1^m$, $m\geq 1$.
 \item[(ii)] Each $V_m^1$ is finite dimensional, $<V_m^1, sL_0^{m+2}>= 0$, and the induced linear maps
 $$s\left( L_0/L_0^{m+1}\right) \stackrel{\cong}{\rightarrow} (V_m^1)^\#$$
 are isomorphisms.
 \item[(iii)] $L_0= \varprojlim_m L_0/L_0^m$ is a pronilpotent Lie algebra.
 \end{enumerate}
 \end{theoremb}

 Next observe that the isomorphism $\pi_*\vert \land V,d\vert \stackrel{\cong}{\rightarrow} \pi_*(\land V,d)$ transports the action of $\pi_1\vert \land V,d\vert$ in $\pi_k\vert \land V,d\vert$ to a representation of $\pi_1 $ in $\pi_k(\land V,d)$. Denote this by
 $$\alpha \times \beta \mapsto \alpha\bullet \beta\,, \hspace{1cm} \alpha \in \pi_1, \beta\in \pi_k(\land V,d)\,.$$
 
 On the other hand,  suppose $x\in L_0$,  $y\in L_k$ and  $k\geq 1$. Then it follows from the Sullivan condition  that for some $m=m(y)$, $(ad\, x)^m(y)= 0$. Thus, when dim$\, H^1(\land V,d)<\infty$, a representation of $G_L$ in $L_k$ is given by 
 $$\exp (x)\bullet y= e^{ad\, x}(y)\,.$$
 Moreover, in this case,
 $$\alpha\bullet \beta = s\,\left(e^{ad\, (s^{-1}\alpha)}(s^{-1}\beta)\right)\,, \hspace{1cm} \alpha \in \pi_1\,, \beta\in \pi_k(\land V,d)\,.$$

\vspace{3mm} Finally, associated with any differential graded Lie algebra $(E, \partial)$,  $E= E_{\geq 0}$, is the chain coalgebra $C_*(E, \partial) := (\land sE, \partial +\delta)$ whose differential is determined by the conditions $\partial sx= -s\partial x$ and $\delta (sx\land sy) = (-1)^{deg\, ,x+1} s[x,y]$. The dual, $C^*(E, \partial) := [C_*(E, \partial)]^\#$ is a cdga. Moreover, if $E/[E,E]$ is a graded vector space of finite type then so is each $E^n/E^{n+1}$, and in this case
$$(\land V_E,d) = \varinjlim_m C^*(E/E^m)$$
is a quadratic Sullivan algebra. If also $\partial = 0$ then 
$$L_{V_E} = \widehat{E} := \varprojlim_m E/E^m\,.$$

 \vspace{3mm} In particular, in \cite{Ge} Getzler  defines a realization functor for differential graded Lie algebras   that is directly related to   spatial realization. If $(\land V,d)$ is a minimal Sullivan algebra where $V = V^{\geq 2}$ is a graded vector space of finite type, then there is   a dgl $(L,\partial)$ in which $L= L_{\geq 1}$ is a vector space of finite type,     the cochain algebra  $C^*(L,\partial)$ is itself a Sullivan algebra, and $(\land V,d)$ is its minimal Sullivan model. Thus $\vert \land V,d\vert \simeq \vert C^*(L,d)\vert$.

Now we can form the simplicial dgl $L\otimes A_{PL}$, where the bracket is defined by $$[\ell \otimes a, \ell'\otimes a']= (-1)^{deg\, \ell'\cdot deg\, a} [\ell, \ell']\otimes aa'\,.$$

Now, for a dgl  $(L,d)$ denote by MC$(L,d)$  the space of Maurer-Cartan elements of $L$, i.e., elements $\omega\in L$ with $d\omega + \frac{1}{2}[\omega, \omega]= 0$.

\begin{theoremb} \label{Geth} When $A$ is an artinian finite dimensional algebra, the natural isomorphism between $\mbox{Hom}
((sL)^\#, A )$ 
and $L\otimes A $ induces a bijection between the space of morphisms from $C^*(L)$ to $A$ and the space MC$(L\otimes A)$. \\
In particular there is an isomorphism of simplicial sets
$$<C^*(L,\partial)> \cong MC(L\otimes A^*_{PL})\,.$$
\end{theoremb}

  \vspace{3mm} Minimal Sullivan models $(\land V,d)$ provide an interesting example of discrete groups, about which a number of open questions remain, including:
  \begin{enumerate}
  \item[1.] When dim$\, H^1(\land V,d)<\infty$ does the inclusion $G_L\hookrightarrow \widehat{UL_0}$ extend to an isomorphism $\widehat{\mathbb Q[G_L]} \stackrel{\cong}{\longrightarrow} \widehat{UL_0}$ from the completion of the group ring.
  \item[2.] Must it always be the case that some $L_0^n /[L_0^n, L_0^n]$ is infinite dimensional ?
  \end{enumerate}
 
 \vspace{3mm}\noindent {\bf Example.} Let $(\land V,d)$ be the minimal Sullivan model of $S^1\vee S^1$. Then $L_V$ is the completion of the free Lie algebra on two generators $E=\mathbb L(x,y)$: $$L_V = \varprojlim_n E/E^n\,.$$
 The group $\pi_1(\land V,d)$ is the associated group and can be described as the vector space $L_V$ equipped with the Baker-Campbell-Hausdorff product,
 $$a*b = a+b + \frac{1}{2} [a,b] + \frac{1}{12} [a,[a,b]] - \frac{1}{12}[b, [a,b]] + \cdots $$

  \section{Whitehead products and nilpotence}
  
  Fix a minimal Sullivan algebra $(\land V,d)$ with homotopy Lie algebra $L$ and, for simplicity, denote $\pi_k(\land V,d)$ simply by $\pi_k$. Recall that the Whitehead products in $\vert \land V,d\vert$ are maps
  $$[\hspace{2mm},\hspace{2mm}]_W : \pi_k(\vert\land V,d\vert) \times \pi_\ell (\vert \land V,d\vert) \to \pi_{k+\ell -1}(\vert\land V,d\vert)\,.$$
  The isomorphism $\pi_k(\vert\land V,d\vert) \to \pi_k(\land V,d)$ then transports these to Whitehead products $[\hspace{2mm},\hspace{2mm}]$ in $\pi_k$. In particular, when $k, \ell \geq 1$,
  $$[sx, sy]_W = (-1)^{deg\, x} s[x,y]\,, \hspace{1cm} x\in L_{k-1}, y\in L_{\ell -1}\,.$$
  
  When $k= 1$ and $\ell = 1$, then
  $$[\alpha, \beta]_W = \alpha\beta\alpha^{-1}\beta^{-1}\,, \hspace{1cm} \alpha, \beta\in \pi_1\,,$$
  and when $k=1$ and $\ell>1$ then
  $$[\alpha, \beta]_W = \alpha\bullet\beta - \beta\,, \hspace{1cm} \alpha\in \pi_1, \beta\in \pi_k\,.$$
  Thus if dim$\, H^1(\land V,d)<\infty$ these can be expressed directly via the exponential map $L\stackrel{\cong}{\rightarrow} G_L$ using the formulae in Section 4.
  
  The three filtrations in Section 4 also generalize to general $k$. First, the lower central series of $\pi_k$ is the sequence $\pi_k = \pi_k^1\supset \cdots \supset \pi_k^r\supset \cdots$, where $\pi_k^{r+1}$ is the (abelian for $k\geq 2$) subgroup generated by the elements $[\alpha, \beta]_W$, $\alpha\in \pi_1, \beta \in \pi_k^r$. Then, the lower central series for $L_k$ is the sequence $L_k = L_k^1 \supset \cdots \supset L_k^r\supset \cdots$ where $L_k^{r+1}$ is the linear span of the Lie brackets $[\alpha, \beta]$, $\alpha \in L_0$, $\beta\in L_k^r$. Finally, denote by $\delta$ the derivation of $\land V$ in which $\delta v$   is the component of $dv$ in $V^1\land V$.  Then a filtration $V^k_{0}\subset V_{1}^k \subset \cdots \subset V^k_{r}\subset \cdots $ of each $V^k$ is given by $V^k_{0} = V^k\cap \mbox{ker}\, \delta$ and $V_{r+1} = \delta^{-1}(V^1\cap V_{r}^k)$. Note that when $k=1$ these filtrations coincide with those defined in Section 4.
  
  \vspace{3mm}\noindent {\bf Definition.} Let $(\land V,d)$ be a minimal Sullivan algebra. Then, with the notation above, nil$\, \pi_k$ ($k\geq 1$), nil$\, L_k$ ($k\geq 0$), and nil$\, V^k$ ($k\geq 1$) are respectively the greatest integers $N_\pi$, $N_L$ and $N_V$ (or $\infty$) such that 
  $$\pi_k^{N_\pi}\neq 0\,, \hspace{5mm} L_k^{N_L}\neq 0\,, \hspace{5mm}\mbox{and } V^k/V^k_{N_V-1}\neq 0\,.$$

  \vspace{3mm} When dim$\, H^1(\land V,d)<\infty$, Section 4 establishes close relations for the filtrations of $\pi_1$, $L_0$ and $V^1$. Similar relations hold for the filtrations of $\pi_k$, $L_{k-1}$ and $V^k$ for an arbitrary $k$, provided   also that dim$\, V^k<\infty$. In particular, these relations hold for all $k$, if dim$\, V<\infty$. Since $(\land V,d)= \varinjlim_\alpha (\land V_\alpha,d)$ with dim$\, V_\alpha<\infty$ an (inverse + direct) limit establishes
  
  \begin{theoremb} \label{nil}For any minimal Sullivan algebra $(\land V,d)$ and any $k\geq 1$,
  $$\mbox{nil}\, \pi_k = \mbox{nil}\, L_{k-1}= \mbox{nil}\, V^k\,.$$
  \end{theoremb}

  \vspace{3mm}\noindent {\bf Remark.} A connected CW complex $X$ is \emph{nilpotent} if $\pi_1(X)$ and each $\pi_k(X)$, $k\geq 2$, are respectively a nilpotent group and a nilpotent $\pi_1(X)$-module. Thus Theorem \ref{nil} characterizes the possible nilpotence of $\vert \land V,d\vert$:
  $$\vert \land V,d\vert \mbox{ is nilpotent }\Longleftrightarrow \mbox{for each $k\geq 1$ and some $N(k)$, $V^k = V^k_{N(k)}$}\,.$$

 \section{Hurewicz homomorphisms}
  
 The Hurewicz homomorphism in topology
 $$ hur : \pi_*(X)\to H_*(X)$$ also has an analogue for minimal Sullivan algebras, $(\land V,d)$. The surjection $\zeta : \land^+V\to \land^+V/\land^+V{\scriptstyle \land} 
  \land^+V = V$   induces a linear map $H(\zeta) : H^+(\land V,d)\to V$, and    for any based CW complex $X$, and a morphism $\varphi : (\land V,d)\to A_{PL}(X)$, the diagram
  $$\xymatrix{
  \pi_*(X) \ar[dd]^{hur}\ar[rrr]^{\pi_*(\varphi)} &&& \pi_*(\land V,d)\ar[dd]^{   H(\zeta)^\#}\\
  \mbox{} \\ H_*(X) \ar[r] & H(X)^\# \ar[rr]^-{ H(\varphi)^\#} && H(\land V,d)^\#}$$
  commutes.
  
 For $n\geq 2$  this is in (\cite{FHTII}). When $n = 1$,   let $\varphi_0$ be the restriction of $\varphi$ to $ V_0^1= V^1\cap$ ker$\, d$.  Then
  $$H(\zeta)^\#\circ \pi_1(\varphi) = \pi_1(\varphi_0) : \pi_1(X) \to \pi_1(\land V_0^1) = (V_0^1)^\#\,.$$
  The discussion of Whitehead products  shows that $\pi_1(\land V_0^1)$ is abelian. Thus $\pi_1(\varphi_0)$ factors through $\pi_1(X)/[\pi_1(X), \pi_1(X)]= H_1(X)$, and so the diagram above commutes  when $n= 1$. 
  
  \vspace{3mm}\noindent {\bf Example 1.}
 If $V = V^k$ then $$\pi_m\vert \land V,d\vert = \left\{ \begin{array}{ll} V^\#\,,  & m= k\\
  0\,, & \mbox{otherwise}.\end{array}\right.$$
  Thus $\vert \land V,d\vert = K(V^\#, k)$.

\vspace{3mm}\noindent {\bf Example 2.}  If $\varphi : (\land V,d)\to A_{PL}(X)$ is the minimal Sullivan model of an $n$-connected CW complex with finite Betti numbers, then $V = V^{\geq n+1}$ is a graded vector space of finite type and 
  $$\pi_k(\varphi) : \pi_k(X)\otimes \mathbb Q \stackrel{\cong}{\longrightarrow} V^\#\,.$$
  
  Since $d: V\to \land^{\geq 2}V$, it vanishes in $V^{\leq 2n}$ and so $H(\zeta) : H^k(\land V)\stackrel{\cong}{\to} V^k$, $k\leq 2n$. Thus the Hurewicz diagram gives
  $$hur : \pi_k(X)\otimes \mathbb Q \stackrel{\cong}{\longrightarrow} H_k(X;\mathbb Q)\,, \hspace{1cm} k\leq 2n\,.$$
  Moreover, in degrees $k\leq 3n+1$, $d : V^k\to \land^2V$. It follows that Im$\, H(\zeta)^k$ is the subspace of $V^k$ orthogonal to the linear span of Whitehead products. This then implies that
  $$\mbox{ker}(hur) : \pi_k(X)\otimes \mathbb Q \to H_*(X;\mathbb Q)\,, \hspace{1cm} k\leq 3n+1$$
  is the linear span of the Whitehead products.

\section{Rationalizations, Malcev completions and Sullivan spaces}

With a nilpotent connected CW complex, $X$, is associated (\cite{FHT}) a continuous map $f : X\to X_{\mathbb Q}$ in which $X_{\mathbb Q}$ is a nilpotent CW complex,  and each of $H_i(f;\mathbb Z)$, $i\geq 1$ and $\pi_i(f)$, $i\geq 2$ extend to isomorphisms
$$
H_i(X;\mathbb Z)\otimes \mathbb Q \stackrel{\cong}{\longrightarrow} H_i(X_{\mathbb Q}, \mathbb Z) \hspace{1cm}\mbox{ and } \pi_i(X)\otimes \mathbb Q \stackrel{\cong}{\longrightarrow} \pi_i(X_{\mathbb Q}).
$$
The space $X_{\mathbb Q}$ is the \emph{rationalization} of $X$, and $X$ is a \emph{rational CW complex} if $f$ is a homotopy equivalence. In particular, a rationalization induces an isomorphism of minimal Sullivan models.

If $\varphi : (\land V,d)\to A_{PL}(X)$ is a minimal Sullivan model of a nilpotent connected CW complex, and in addition $H(X)$  is a graded vector space of finite type, then $\widehat{\varphi} : X\to \vert\land V,d\vert$ is a rationalization, and so $\pi_i(\varphi) : \pi_i(X)\otimes \mathbb Q\stackrel{\cong}{\longrightarrow} V_i^\#$, $i\geq 2$. 

\vspace{3mm}\noindent {\bf Example: Eilenberg-MacLane spaces.} If $G$ is an abelian group and $r\geq 1$, then the Eilenberg-MacLane space $K(G,r)$ is an $H$-space and therefore nilpotent. Write $G_{\mathbb Q}= G\otimes_{\mathbb Z}\mathbb Q$. Then
$$K(G,r)\to K(G_{\mathbb Q}, r)$$
is a rationalization.

Now suppose dim$\, G_{\mathbb Q}<\infty$. Then the minimal Sullivan model of $K(G,r)$ has the form $(\land V^r,0)$ and dim$\, V^r =$ dim$\, G_{\mathbb Q}$. In this case the map
$$K(G_{\mathbb Q},r)\to \vert \land V^r\vert$$
is a homotopy equivalence.

\vspace{3mm} \emph{Sullivan spaces} are a broader class of connected CW complexes $X$ in which the fundamental group, $G$ may not be nilpotent, but in which if $\varphi : (\land V,d)\to A_{PL}(X)$ is a minimal Sullivan model then $\pi_1(\widehat{\varphi}) $ is a \emph{Malcev completion} (described next) and $\pi_i(\widehat{\varphi})$, $i\geq 2$ and $H^i(\widehat{\varphi})$, $i\geq 0$ are isomorphisms,
$$\pi_i(X)\otimes \mathbb Q \stackrel{\cong}{\longrightarrow} \pi_i\vert\land V,d\vert \hspace{1cm}\mbox{ and } H^i(\vert\land V,d\vert) \stackrel{\cong}{\longrightarrow} H^i(X)\,.$$

Recall then $G^2 = [G,G]$ is the commutator subgroup of $G$. 
Because $H^1(X)= \left( G/[G,G]\right)^\#$ it follows from the Hurewicz diagram that
$
\mbox{dim}\, H^1(X)<\infty $ if and only if the morphism induced by $\pi_1(\varphi)$ is an isomorphism.
$$ G/[G,G]\otimes \mathbb Q \stackrel{\cong}{\to} \frac{\pi_1(\land V,d)}{[\pi_1(\land V,d), \pi_1(\land V,d)]}\,.$$

When this holds $\pi_1(\varphi)$ induces isomorphisms
$$G^n/G^{n+1}\otimes \mathbb Q \stackrel{\cong}{\longrightarrow} \pi_1^n(\vert \land V,d\vert)/ \pi_1^n(\vert \land V,d\vert)\,, \hspace{3mm} n\geq 1\,,$$
and  
$$\pi_1(\land V,d) \stackrel{\cong}{\longrightarrow } \varprojlim_n \pi_1(\land V,d)/\pi_1^n(\land V,d)\,.$$
These isomorphisms identify $\pi_1(\land V,d)$ as a \emph{Malcev completion} of $G$. 

\vspace{3mm} The distinct role of fundamental groups in topology as opposed to the abelian higher homotopy groups is reflected in Sullivan's models: if $\varphi : (\land V,d)\stackrel{\simeq}{\to} A_{PL}(X)$ is a minimal Sullivan model then while $H(\varphi)$ is an isomorphism $\pi_*(\varphi)$ may not be. Now, the fundamental group $G$ of $X$ determines a classifying map $\gamma : X\to K(G,1)= BG$, whose homotopy fibre is the universal covering space $\widetilde{X}$ of $X$. Then, because $H^1(\gamma)$ and $H^2(\gamma)$ are respectively an isomorphism and injective, there is a commutative diagram,
$$\xymatrix{
(\land V^1,d) \ar[d]_{\varphi_G}\ar[rr]^\lambda && (\land V,d)\ar[d]_\varphi^\simeq\ar[rr]^\rho && (\land V^{\geq 2}, \overline{d})\ar[d]^{\overline{\varphi}}\\
A_{PL}(BG) \ar[rr]^{A_{PL}(\gamma)} && A_{PL}(X)\ar[rr] && A_{PL}(\widetilde{X})\,,}$$
in which $\varphi$ is the minimal Sullivan model of $X$. Moreover, $(\land V^1,d)$ and the morphism $\varphi_G$ are invariants just of $G$. 

Thus for $\varphi$ to fully capture the "rational homotopy type" of $X$ it should satisfy the following three conditions:
\begin{enumerate}
\item[$\bullet$] $\varphi_G$ is a quasi-isomorphism
\item[$\bullet$] $\pi_1(\varphi)$ is a Malcev completion
\item[$\bullet$] For $n\geq 2$, $\pi_n(\varphi) : \pi_n(X)\otimes \mathbb Q \stackrel{\cong}{\to} \pi_n(\land V,d)$.
\end{enumerate}

\vspace{3mm}\noindent {\bf Definition.} $X$ is a \emph{Sullivan space} if the three conditions above are satisfied by its minimal Sullivan model.

\vspace{3mm}\noindent {\bf Remark.} The second condition   implies that dim$\, H^1(X)<\infty$ and it follows that $V^1$ is at most countable. If $V^{<n}$ is at most countable then the Hurewicz diagram together with the third condition   implies that $V^n$ is at most countable via the exact sequence
$$V^n\cap \mbox{ker}\, d \to V\stackrel{d}{\to} \land V^{<n}\,.$$
Thus these two conditions imply that $\land V$ is countable and hence that $H(X)$ is countable. But a cohomology algebra cannot be countable unless it has finite type. Thus
$$X \mbox{ is a Sullivan space} \Longrightarrow H(X) \mbox{ has finite type}.$$

It follows that  this definition is equivalent to the definition in \cite{FHTII}, and so $X$ is a Sullivan space if and only if : 
\begin{enumerate}
\item[$\bullet$] $H(X)$ is a graded vector space of finite type;
\item[$\bullet$] $H(BG) = \varinjlim_n H(B\, G/G^n)$;
\item[$\bullet$] $G$ acts nilpotently by covering transformations in $H(\widetilde{X})$.
\end{enumerate}

\vspace{3mm}
Finally, the minimal Sullivan model $(\land V,d)$ of a Sullivan space satisfies
\begin{enumerate}
\item[$\bullet$] $H(\land V,d)$ and $V^{\geq 2}$ are graded vector spaces of finite type,
\item[$\bullet$] $H(BG)= \varinjlim_n H(B\, G/G^n)$, $G$ denoting $\pi_1(\land V,d)$.
\end{enumerate}
Conversely, if a minimal Sullivan model satisfies these two conditions, then
$$m_V : (\land V,d) \to A_{PL}(\vert \land V,d\vert)$$
is a quasi-isomorphism and $\vert \land V,d\vert$ is a Sullivan space. In particular,

\begin{theoremb} If $X$ is a Sullivan space with minimal Sullivan model $(\land V,d)$, then 
$H(\widehat{\varphi})$ is an isomorphism, $\pi_1(\widehat{\varphi})$ is a Malcev completion, and  $$\pi_i(\widehat{\varphi}) : \pi_i(X)\otimes \mathbb Q \stackrel{\cong}{\longrightarrow} \pi_i(\vert\land V,d\vert)\,, \hspace{1cm} i\geq 2\,.$$
\end{theoremb}

\vspace{3mm}\noindent {\bf Remark.} The spatial realization $\vert \land V,d\vert$ of a minimal Sullivan algebra is a nilpotent Sullivan space if and only if dim$\, V^p<\infty$, $p\geq 1$. In this case, $H_{>0}(\vert \land V,d\vert; \mathbb Z)$ is a rational vector space.

It is an open question as to when in general the realization $\vert \land V,d\vert$ satisfies $H_{>0}(\vert \land V,d\vert;\mathbb Z)$ is a rational vector space even in the case $(\land V,d)$ is the minimal Sullivan model of a Sullivan space.

\vspace{3mm}\noindent {\bf Examples.} 
\begin{enumerate}
\item[1.] A  simply connected space $X$ is a Sullivan space if and only if $H(X)$ has finite type.
\item[2.] Nilpotent connected CW complexes with finite Betti numbers are Sullivan spaces.
\item[3.] An oriented Riemann surface is a Sullivan space.
\item[4.] A finite wedge of circles is a Sullivan space.
\item[5.] The classifying spaces of pure braid groups and the Heisenberg group are Sullivan spaces.
\item[6.] The connected sums $\mathbb RP^{2k+1}\#\mathbb RP^{2k+1}$, $k>0$, are not Sullivan spaces, but the classifying spaces of their fundamental groups are Sullivan spaces.
\end{enumerate}

\vspace{3mm} It is an open question whether $\vert \land V,d\vert$ is a Sullivan space if $m_V : (\land V,d)\to A_{PL}\vert \land V,d\vert$ is a quasi-isomorphism.

 \section{LS category}

 The LS (Lusternik-Schnirelmann) category of a topological space $X$ is the least $m$ (or $\infty$) for which $X$ can be covered by $m+1$ open sets, each contractible in $X$. If $X$ is a normal space then cat$\, X$ is the least $m$ for which $X$ is a retract of an $m$-cone.
 
 Analogously, a Sullivan algebra $(\land V,d)$ determines diagrams,
 $$\xymatrix{(\land V,d)\ar[d]^\rho\\
 (\land V/\land^{>m}V, \overline{d}) && (\land W,d)\ar[ll]^{\simeq}_\varphi\ar@{-->}[llu]_\psi\,, & m\geq 0,}$$
 in which $\rho $ is the quotient map and $\varphi$  is a Sullivan model. The \emph{LS category, cat$\, (\land V,d)$} of $(\land V,d)$ is then the least $m$ (or $\infty$) such that there is a morphism $\psi : (\land W,d)\to (\land V,d)$ satisfying $\rho\circ \psi \sim \varphi$. 
 
 This invariant has three key properties, which together are central to a number of applications in the following sections of this survey. Thus, let $(\land V,d)$ be any minimal Sullivan algebra.  
 \begin{enumerate}
 \item[(i)] If $(\land V,d)$ is the minimal Sullivan model of a connected CW complex $X$ then $$\mbox{cat}\, (\land V,d) \leq \mbox{cat}\, X\,.$$
 \item[(ii)] (Mapping Theorem). If $(\land V,d)\to (\land Z,d)$ is a surjection of minimal Sullivan algebras then
 $$\mbox{cat}\, (\land V,d) \geq \mbox{cat}\, (\land Z,d)\,.$$
 \item[(iii)] If $H^{>n}(\land V,d) = 0$, then cat$\, (\land V,d)\leq n$.
 \end{enumerate}
 
\vspace{3mm} Now fix a minimal Sullivan algebra $(\land V,d)$. The invariants mcat$(M,d)$ and $e(M,d)$ are then defined   for any $(\land V,d)$-module $(M,d)$, via a construction analogous to the definition of cat$\, (\land V,d)$   : associated with $(M,d)$ are diagrams of $(\land V,d)$-modules
 $$\xymatrix{(P,d)\ar[d]^\rho\\
 (P/\land^{>m}V,\overline{d}) && (Q,d)\ar[ll]^\simeq_\varphi}$$
 in which $(P,d)$ is a semi-free model of $(M,d)$, $\rho$ is the quotient map and $\varphi$ is a semifree model. 
  The \emph{module category, mcat$(M,d)$}  is then the least $m$ (or $\infty$) such that there is a morphism of modules $\psi : (Q,d)\to (P,d)$ such that $\rho\circ \psi\sim \varphi$, and  the \emph{Toomer invariant} $e(M,d)$ is the least $m$ (or $\infty$) for which $H(\rho)$ is injective. They satisfy
  $$e(\land V,d) \leq \mbox{mcat}\, (\land V,d) = \mbox{cat}\, (\land V,d)\,,$$
  the  inequality following because 
  $(\land V,d)$ is $(\land V,d)$-semifree,   while the equality is a deep theorem of Hess (\cite{Hess}).
 
The Toomer invariant can also be defined for a connected CW complex $X$, and an  old question of Berstein asked whether rationally $e(X)= \mbox{cat}\, X$. This turns out not to be true, but it   follows from Hess' theorem that if   $H(\land V,d)$ is a vector space of finite type then
 $$\mbox{cat}\, (\land V,d) = e((\land V,d)^\#)\,.$$
 
 This has two consequences:
 \begin{enumerate}
 \item[$\bullet$] If both $(\land V,d)$ and $(\land W,d)$ are minimal Sullivan algebras with homology of finite type, then
 $$\mbox{cat}\, ((\land V,d)\otimes (\land W,d)) = \mbox{cat}\, (\land V,d) + \mbox{cat}\, (\land W,d)\,.$$
 \item[$\bullet$] If $H(\land V,d)$ satisfies Poincar\'e duality then
 $$e(\land V,d) =\mbox{cat}\, (\land V,d)\,.$$
 \end{enumerate}
 
 \vspace{3mm}\noindent {\bf Remark 1.} If $X$ is a simply connected CW complex with rational homology groups of finite type and minimal Sullivan model $(\land V,d)$ then cat$\, X_{\mathbb Q}=$ cat$(\land V,d)$. As pointed out to us by Parent, the hypothesis of simple connectivity is essential: if $\land V= \land v$, deg$\, v= 1$, then cat $(\land V)= 1$ but cat$\, \vert \land V\vert = 2$. 
 
 {\bf 2.} The theory of semi-free resolutions, and with it the definition of the Toomer invariant, applies to modules over any dga. In particular, suppose $(\land V,d)$ is the minimal Sullivan model of a connected CW complex $X$ with finite Betti numbers.
 Then, as for any path connected space, $(\land V,d)$ is connected by quasi-isomorphism to the singular cochain complex $C^*(X;\mathbb Q)$. It follows that
 $$  e(C_*(X;\mathbb Q))= \mbox{cat}\, (\land V,d)\,,$$
 where $C^*(X;\mathbb Q)$ acts on the  chain complex $C_*(X;\mathbb Q)$ by the cap product.
 
 \section{The trichotomy : growth of homotopy groups}
 
  For any based  connected CW complex, set
$$\mbox{rk}_k(X) = \mbox{dim}\, \pi_k(X)\otimes \mathbb Q\,, \hspace{1cm} k\geq 2\,.$$
There are then exactly three mutually exclusive possibilities:
\begin{enumerate}
\item[$\bullet$] $\sum_{k\geq 2} \mbox{rk}_k(X)<\infty$  ($X$ is rationally elliptic).
\item[$\bullet$] Each $\mbox{rk}_k(X)<\infty$ but $\sum_{k\geq 2}\, \mbox{rk}_k(X)=\infty $ ($X$ is rationally hyperbolic).
\item[$\bullet$] Some $\mbox{rk}_k(X)= \infty$ ($X$ is rationally infinite).
\end{enumerate}

\vspace{3mm} When $X$ is  $N$-dimensional an application of Gottlieb's evaluation subgroups (\cite{Go}) establishes in the rationally infinite case that there is a constant $K$ and an infinite sequence $\cdots k_n<k_{n+1}<\cdots $ such that $k_n\leq nK$ and $\mbox{rk}_{k_n}(X) = \infty$.  

When $X$ is rationally elliptic or rationally hyperbolic then the universal cover $\widetilde{X}$ is a Sullivan space and its minimal model $(\land V,d)$ satisfies
$$\mbox{dim}\, V^k = \mbox{rk}_k(X)\,.$$
Moreover, since cat$\, \widetilde{X}\leq \mbox{cat}\, X$ and dim$\, \widetilde{X}=$ dim$\, X$,  it follows that if cat$\, X<\infty$ (resp. dim$\, X\leq N$) then cat$\, (\land V,d)<\infty$ (resp. $H^{> N}(\land V,d) = 0$). 
In these two cases a non-elementary application of minimal Sullivan model technology yields strong information on the size and growth of the integers rk$_k(X)$. 

\vspace{3mm}\noindent {\bf 1. }{\sl The rationally hyperbolic case.}\\

\vspace{2mm} When $X$ is rationally hyperbolic,
$$\alpha_X := \lim\sup \frac{\log \mbox{ rk}_k(X)}{k} >0\,,$$
and the $\mbox{rk}_k(X)$ exhibit the following  very refined version of exponential growth: for some positive integer $r$
$$\varinjlim_{n\to \infty} \max_{n\leq k\leq n+r} \frac{\log \mbox{ rk}_k(X)}{k}  = \alpha_X\,.$$
When $H^{>N}(\widetilde{X}) = 0$,   there are explicit error terms for the limit above, and $r$ can be taken as $N-2$:  there are positive constants $\gamma, \beta$ depending only on  $H(\widetilde{X})$,  and there is a positive constant $K$ such that
$$\alpha_X - \frac{\gamma}{\log n} \leq \max_{n\leq k\leq n+N-2} \frac{\log \mbox{rk}_k(X)}{k}  \leq \alpha_X + \frac{\beta}{n}\,, \hspace{5mm} n\geq K\,.$$
(Note the typo in the statements of Theorems 12.10 and 12.11 in \cite{FHTII} where the left hand error term incorrectly appears as $\gamma \frac{\log n}{n}$ instead of $\frac{\gamma}{\log n}$.)

\vspace{3mm}
For a family of rationally hyperbolic spaces containing spaces of LS category 1 and 2, a stronger estimate of $\mbox{rk}_k(X)$ is given by Lambrechts \cite{La}: There are constants $A$ and $B$ such that 
$$\frac{A}{n} e^{n\alpha_X} \leq \max_{n\leq k\leq n+N-2}   \mbox{rk}_k(X) \leq \frac{B}{n}e^{n\alpha_X}\,.$$  

There are open questions about rationally hyperbolic spaces:
\begin{enumerate}
\item[1.] Does the stronger estimate of Lambrechts hold in general ?
\item[2.] Do the numbers dim$\, (L_V/[L_V, L_V])^k$ grow at most exponentially at a slower rate then the numbers rk$_k(X)$ ?
\item[3.] There is an analogous phenomenon in the theory of local commutative noetherian rings $R$, where deviations $\varepsilon_k$ play the role of the numbers rk$_k(X)$. If $R$ is not a local complete intersection it is known that each $\varepsilon_k\neq 0$, but whether they satisfy the same growth as the rk$_k(X)$ in the hyperbolic case is unknown.
\end{enumerate}

\vspace{2mm}
\noindent {\bf 2.}  {\sl The rationally elliptic case.} 

\vspace{2mm} When $X$ is rationally elliptic,  if cat$\, X<\infty$ then $H(\widetilde{X})$ is finite dimensional, so that for some $N$, $H^{>N}(\widetilde{X})= 0$. In this case $\sum_i \mbox{dim} H_i(\widetilde{X}) \leq 2^N$, and 
$$r_k(X) = 0,  \hspace{3mm} k\geq 2N \hspace{3mm}\mbox{and } \sum_{k\geq 2} r_k(X) \leq 2 \mbox{ cat}\, X\,.$$
Moreover, there is an explicit algorithm (\cite{FH}) which decides when a sequence $r_2, \cdots , r_{2N-1}$ arises as the ranks of the rational homotopy groups of a rationally elliptic spaces $X$. In fact a finite sequence $2a_1, \cdots , 2a_q, 2b_1-1, \cdots , 2b_p-1$ are the degrees of a basis of $\pi_{\geq 2}(X)\otimes \mathbb Q$ with $X$ rationally elliptic and  cat$\, X<\infty$ if and only if the following condition  is satisfied:
\begin{enumerate}\item[$\bullet$] For any subsequence $a_{j_1}\cdots a_{j_s}$ at least $s$ of the $b_i$ can be written as an integral combination $b_i = \sum k_\lambda a_{j_\lambda}$ with each $k_{\lambda}\geq 0$ and $\sum k_\lambda \geq 2$.\end{enumerate}

\vspace{2mm} If $X$ is a nilpotent, rationally elliptic space then its minimal Sullivan model $(\land V,d)$ satisfies dim$\, V<\infty$. If in addition, cat$\, X<\infty$ then $H(X)$ is a finite dimensional Poincar\'e duality algebra. In this case, the homotopy Euler characteristic $\chi_\pi := \mbox{dim}\, V^{even}- \mbox{dim}\, V^{odd}$ satisfies $\chi_\pi \leq 0$ and 
$$\chi_\pi = 0\Longleftrightarrow H^{odd}(X) = 0\,.$$
Further, if $\chi_\pi <0$ then the Euler characteristic of $H(X)$ is zero. Finally, if $\chi_\pi = 0$, then $X$ is a formal space and $(\land V,d)$ has the form $(\land (v_1, \cdots , v_n, w_1, \cdots , w_n),d)$ with deg$\, v_i$ even, deg$\, w_j$ odd, $dv_i= 0$ and $dw_j\in \land (v_1, \cdots , v_n)$. In this case $$H(X) = \land (v_1, \cdots , v_n)/(dw_1, \cdots , dw_n)\,.$$

Nilpotent  connected rationally elliptic spaces $X$ with cat$\, X<\infty$ include compact connected  Lie groups $G$, homogeneous spaces $G/K$ (with $G, K$ compact and connected), and connected sums $\mathbb CP^n \# \mathbb CP^n$. Finally if a compact connected Lie group $G$ acts on a closed manifold $M$ with principal orbits of codimension 1, or if $M$ is a closed Dupin hypersurface in $S^{n+1}$, then $M$ is rationally elliptic (\cite{GH}).

There remain a number of open questions about rationally elliptic spaces:
\begin{enumerate}
\item[1] In \cite{Bott} Berger and Bott   limit the rate of potential exponential growth of the Betti numbers of a closed positively curved manifold. An open conjecture, attributed to Bott, asserts that such manifolds are rationally elliptic.
\item[2.] If $X$ is rationally elliptic, cat$\, X<\infty$ and $\chi_\pi = 0$ it remains an open conjecture whether Der$_{<0}(H(X))= 0$. This has been resolved when $X$ is a homogeneous space (\cite{Shiga}).
\item[3.] An open conjecture of Hilali (\cite{Hi}) asserts that if $X$ is rationally elliptic and cat$\, X<\infty$, then $$\sum_k \mbox{rk}_k(X) \sum_k \mbox{dim}\, H^k(X)\,.$$
 \end{enumerate}

\section{Minimal $\land$-models}

Any continuous map $f : X\to Y$ may be converted to a fibration in a commutative diagram of the form
$$\xymatrix{X' \ar[rr] \ar[d]_p && X\ar[dll]^f\\
Y}$$
in which $p$ is a fibration and the horizontal arrow is a homotopy equivalence. In this case the fibre of $p$ is a homotopy invariant of $f$, called its \emph{homotopy fibre}. 

There is an analogue construction for morphisms $\psi : (A,d)\to (C,d)$ of cdga's satisfying $H^0(B)= \mathbb Q= H^0(C)$.  

\vspace{3mm}\noindent {\bf Definition 1.} A \emph{minimal $\land$-extension} of a cdga $(B,d)$ is a cdga morphism of the form $$\lambda : (B,d)\to (B\otimes \land Z,d)$$
in which $\lambda (b) = b\otimes 1$, $Z = Z^{\geq 0}$, and there is an increasing filtration $0= Z(0) \subset \cdots \subset Z(r)\subset \cdots $ for which $Z = \cup_r Z(r)$ and 
$$d : 1\otimes Z^p(r) \to B\otimes \land (Z^{<p}\oplus Z^{\leq p}(r-1))\,.$$
If $Z = Z^{\geq 1}$ this is a \emph{minimal  Sullivan extension}.

\vspace{2mm} {\bf 2.} A \emph{minimal $\land$-model} (resp. a \emph{minimal Sullivan model}) for a cdga morphism $\psi : (B,d)\to (C,d)$ is a commutative diagram
$$\xymatrix{ (B\otimes \land Z,d)\ar[rr]_\simeq^\varphi && (C,d)\\
(B,d)\ar[u]^\lambda \ar[urr]_\psi}$$
in which $\lambda (b)= b\otimes 1$ and $\varphi$ is a quasi-isomorphism from a minimal $\land$-extension (resp. a minimal Sullivan extension).

\vspace{3mm}\noindent {\bf Remark.}  While $Z$ has a filtration satisfying a condition similar to that of the defining condition for Sullivan algebras, because $H^1(\psi)$ may not be an isomorphism, $Z$ may   have elements of degree $0$, which can invalidate some arguments used for minimal Sullivan algebras $(\land V,d)$.

\vspace{3mm}  The definition of homotopy, the homotopy lifting theorem, and the existence and uniqueness of minimal $\land$-models all generalize from the theory of Sullivan models of a cdga. In particular, two cdga morphisms
$$\varphi_0, \varphi_1 : (B\otimes \land Z,d)\to (C,d)$$
are  \emph{homotopic rel $B$} ($\varphi_0\sim_B\varphi_1$) if there is a morphism $\Phi : (B\otimes \land Z,d)\to \land (t,dt) \otimes (C,d)$ for which $\varphi_0$, $\Phi$ and $\varphi_1$ restrict to the same morphism in $B\otimes 1$ and $(\varepsilon_i\otimes id)\circ \Phi = \varphi_i$; there is also the

\begin{theoremb} If $\psi : (B,d)\to (C,d)$ is any cdga-morphism in which $H^0(B) = \mathbb Q = H^0(C)$, then $\psi$ has a minimal $\land$-model,  and any two are isomorphic by an isomorphism restricting to the identity in $B$.\end{theoremb}

\vspace{3mm}\noindent {\bf A closed model category}. (\cite{BG}, \cite{Hsu}) In \cite{Q} Quillen introduced closed model categories as an abstract framework within which homotopy theoretic techniques apply. In particular, a closed model category, ${\mathcal C}$, is defined as follows:  
 \begin{enumerate}
 \item[$\bullet$] ${\mathcal C}$ is the category of complexes $(A,d)$ in which $A = A_{\geq 0}$ has a commutative product, and $d$ is a derivation, but $A$ is \emph{not} required to have an identity $1\in A_0$.
 \item[$\bullet$] $0$ is both the initial and terminal object.
 \item[$\bullet$] The weak equivalences are the quasi-isomorphisms
 \item[$\bullet$] The fibrations are the surjective morphisms
 \item[$\bullet$] Minimal $\land$-extensions $(A,d)\to (A\otimes \land V,d)$ are cofibrations.
 \end{enumerate}
In this context right homotopy coincides with Sullivan's definition.

\section{Minimal $\land$-models of a fibration}

Suppose that $p : X\to Y$ is a fibration of path connected spaces, and $j: F\to X$ is the inclusion of the fibre at a base point $y_0\in Y$. Then a minimal Sullivan model $\varphi_Y: (\land W,d) \to A_{PL}(Y)$ extends to a minimal $\land$-model
$$\varphi_X : (\land W\otimes \land Z,d) \to A_{PL}(X)$$
of the morphism $A_{PL}(p)\circ \varphi_Y$. Because the composite $F\to X\to Y$ is the constant map it follows that $A_{PL}(j)\circ A_{PL}(p)$ is the augmentation $A_{PL}(Y)\to A_{PL}(y_0)= \mathbb Q$. Thus $A_{PL}(j)\circ \varphi_X$ factors to give the commutative diagram,
$$\xymatrix{
(\land W,d) \ar[rr]^\lambda 
\ar[d]^{\varphi_Y}_\simeq 
&& (\land W\otimes \land Z,d)
 \ar[rr]^{\rho}
\ar[d]^{\varphi_X}_\simeq
 && (\land Z, \overline{d})
  \ar[d]^{\varphi_F}
\\
A_{PL}(Y) 
 \ar[rr]^{A_{PL}(p)} && A_{PL}(X)  \ar[rr]^{A_{PL}(j)} && A_{PL}(F)\,,
}$$
even when $F$ is not path connected and $Z^0\neq 0$. 

\vspace{3mm}\noindent {\bf Definition.} This diagram is \emph{the minimal $\land$-model} of the fibration.

\vspace{3mm} Now the homotopy lifting property of a fibration results in a right "action up to homotopy" of $\pi_1(Y, y_0)$ on $F$, which in particular determines a left \emph{holonomy representation} of $\pi_1(Y, y_0)$ in $H(F)$. Thus, while in general $\varphi_F$ may not be a quasi-isomorphism there is the

\begin{theoremb} \label{12.1} Assume in the diagram above that $F$ is path connected, that one of $H(F)$ and $H(Y)$ is a graded vector space of finite type, and that $\pi_1(Y, y_0)$ acts nilpotently in each $H^k(F)$. Then $H(\varphi_F)$ is also an isomorphism.\end{theoremb}

\vspace{3mm}\noindent {\bf Remark.} Theorem \ref{12.1} also holds when $p$ is a Serre fibration and $F$ has the homotopy type of a CW complex.

\vspace{3mm} Conversely, suppose $(\land W\otimes \land Z,d)$ is a Sullivan extension, and that $\varphi_Y : (\land W,d)\to A_{PL}(Y)$ is a minimal Sullivan model. Then the sequence $\vert \land W,d\vert \leftarrow \vert\land W\otimes \land Z,d\vert \leftarrow \vert\land Z,\overline{d}\vert$ is a Serre fibration. Pulling this back over $\widehat{\varphi}_Y$ gives a Serre fibration
$$\xymatrix{ Y & X\ar[l]_p & \vert\land Z,d\vert\ar[l]}$$
and a commutative diagram
$$\xymatrix{
(\land W,d) \ar[rr] \ar[d]^{\varphi_Y} && (\land W\otimes \land Z,d) \ar[rr]\ar[d]^{\varphi_X} && (\land Z, \overline{d})\ar[d]^{m_Z}\\
A_{PL}(Y)\ar[rr] && A_{PL}(X) \ar[rr] && A_{PL}\vert\land Z,\overline{d}\vert\,.}$$
Moreover, if $m_Z$ is a quasi-isomorphism, so is $\varphi_X$.

\vspace{3mm}\noindent {\bf Example.} In the diagram above suppose that $Y$ is a Sullivan space and that $\pi_1(p)$ and each $\pi_k(p)\otimes \mathbb Q$ are surjective. Then $Z^0= 0$, $(\land W\otimes \land Z,d)$ is the minimal Sullivan model of $X$, and 
$$\mbox{cat}\, (\land Z, \overline{d}) \leq \mbox{cat}\, X\,.$$

In fact, let $S\subset Z$ be the subspace such that $d(1\otimes S) \subset \land^{\geq 2}(W\oplus Z)$ and write $Z = S\oplus T$. Then division by $\land^{\geq 2}(W\oplus Z)$ induces from $d$ a linear injection $d_0 : 1\otimes T\to W$ and division by $1\otimes T$  and $d(1\otimes T)$ is a quasi-isomorphism $\psi : (\land W\otimes \land Z,d)$ onto a minimal Sullivan algebra $(\land V,d)$. 

By construction, $(\land W\otimes \land Z,d)$ is the minimal Sullivan model of $X$ if and only if $T= 0$, or equivalently if the morphism $\overline{\lambda} : (\land W,d) \to (\land V,d)$ induced from $\lambda$ is injective. But since $Y$ is a Sullivan space, for each $w\in W$ there is an $\alpha \in\pi_*(Y)$ satisfying $$<w, \pi_*(\varphi_Y)\alpha> \neq 0\,.$$ Then by hypothesis for some $\beta \in \pi_*(X)$, $<w, \pi_*(\varphi_Y)\pi_*(\rho)\beta> \neq 0$. Since $\pi_*(\varphi_Y)\pi_*(\rho) = \pi_*(\overline{\lambda})\pi_*(\varphi_X)$ this gives
$$0\neq <w, \pi_*(\overline{\lambda})\pi_*(\varphi_X)\beta> = <Q(\overline{\lambda})w, \pi_*(\varphi_X)\beta>\,,$$
where $Q(\overline{\lambda})w$ is the component in $V$ of $\overline{\lambda}w$.

It follows that $\overline{\lambda}$ is injective and $(\land W\otimes \land Z,d)$ is the minimal Sullivan model of $X$. In particular, $Z^0 = 0$. Finally, the Mapping Theorem gives
$$\mbox{cat}\, X\geq \mbox{cat}\, (\land W\otimes \land Z,d) \geq \mbox{cat}\, (\land Z,d)\,.$$

\section{Holonomy representations}

Suppose $\xymatrix{F\ar[r]^j & X\ar[r]^p & Y}$ is a fibration in which $j$ is the inclusion of the fibre at $y_0\in Y$. As described in \S 12, there is then a right homotopy action of $\pi_1(Y, y_0)$ on $F$, which in turn gives a left \emph{holonomy representation} $hol : \pi_1(Y, y_0)\to \mbox{Aut}\, (H(F))$. (In \cite{FHTII} this is denoted $H\, hol$)

There is an analogous construction for $\land$-extensions $(\land W,d)\to (\land W\otimes \land Z,d)\to (\land Z, \overline{d})$ of a minimal Sullivan algebra when dim$\, H^1(\land W,d)<\infty$. Here we obtain a holonomy representation
$hol : \pi_1(\land W,d)\to \mbox{Aut}\, H(\land Z,\overline{d})$ constructed from the differential $d$ in $\land W\otimes \land Z$. Moreover, if
$$\xymatrix{
(\land W,d) \ar[d]^{\varphi_Y} \ar[r] & (\land W\otimes \land Z,d) \ar[d]^{\varphi_X}\ar[r]& (\land Z,\overline{d})\ar[d]^{\varphi_F}\\
A_{PL}(Y) \ar[r] & A_{PL}(X) \ar[r] & A_{PL}(F)}$$
is a commutative diagram then for $\alpha\in \pi_1(Y)$,
$$H(\varphi_F)\circ hol\, (\pi_1(\varphi_Y)\alpha) = hol\, (\alpha) \circ H(\varphi_F)\,.$$

The representation $hol$ of $\pi_1(\land W,d)$ is constructed from the \emph{holonomy representation}, $$\overline{\theta} : L_W \to \mbox{Der}\, H(\land Z)\,,$$ defined as follows.
Filter  $\land W\otimes \land Z$ by the ideals $\land^{\geq k}W\otimes \land Z$ to yield  a spectral sequence whose $E_1$-term is $(\land W\otimes H(\land Z,\overline{d}),d_1)$. Then   $d_1(1\otimes [\Phi]) = \sum_i w_i\otimes \overline{\theta}_i ([\Phi])$, $w_i$ denoting a basis of $W$, and
$$\overline{\theta}(x)([\Phi]) = - \sum <w_i, sx> \overline{\theta}_i ([\Phi])\,, \hspace{5mm} x\in L_W.$$
This then extends to the holonomy representation $\overline{\theta}$ of $UL_W$ in $H(\land Z, \overline{d})$. 

In particular, the restriction of $\overline{\theta}$ to $L_{W_0}$ is locally nilpotent and so extends to $\widehat{UL_{W_0}}$. Restricting this to $\pi_1(\land W,d)$ gives the representation $hol$ of that group.

\section{Fiber squares}

A \emph{homotopy pullback diagram} is a map of Serre fibrations,
$$\xymatrix{
E'\ar[rr] \ar[d]^{p'} && E\ar[d]^p\\
B' \ar[rr]^f && B
}$$
inducing a weak homotopy equivalence in the fibres. (Equivalently, $E'\to E\times_BB'$ is a weak homotopy equivalence).

Suppose given such a homotopy pullback diagram in which $B$ and $B'$ are both path connected, the fibre $F$ has finite Betti numbers, and $\pi_1(B)$ acts nilpotently in $H(F)$. 
 Let $(\land W,d)\to (\land W\otimes \land Z,d)$ be a minimal   Sullivan model for $p$ and let $\varphi : (\land W,d)\to (\land V,d)$ be a Sullivan representative for $f$. Then a   Sullivan model for $p'$ is given by
$$(\land V,d)\to (\land V\otimes \land Z,d')\,,$$
where $d'(v)= (\varphi\otimes id) d(v)$.

\vspace{3mm}\noindent {\bf Example 1: Elementary principal fibrations}

For any abelian group $G$, an \emph{Eilenberg-MacLane space} $K(G,r)$, is a path connected based topological space for which 
$$\pi_i(K(G,r)) = \left\{\begin{array}{ll} G\,, \hspace{1cm}\mbox{} & i= r\\0, & \mbox{otherwise}.\end{array}\right.$$
Now suppose $f : Y\to K(G,r)$ is a map from a path connected space to an Eilenberg-MacLane space with $r\geq 2$. If $q: P\to K(G,r)$ is a fibration with contractible total space $P$, and if 
$$\xymatrix{X\ar[rr] \ar[d]^p && P\ar[d]\\
Y\ar[rr]^f && K(G,r)}$$
is a homotopy pullback diagram, then the fibre $F$ of $p$ is a $K(G, r-1)$ and $\pi_1(Y)$ acts nilpotently in $H(F)$. The fibration $p$ is called an \emph{elementary principal fibration}.

If in addition $G_{\mathbb Q}:= G\otimes \mathbb Q$ is finite dimensional then $H(K(G, r)) = \land V^{r}$ and dim$\, V^{r} =$ dim$\, G_{\mathbb Q}$. In this case there is a commutative diagram
$$\xymatrix{(A_{PL}(Y)\otimes \land Z,d) \ar[rr]^\simeq && A_{PL}(X)\\\mbox{} & A_{PL}(Y)\ar[lu] \ar[ru]_{A_{PL}(p)} \,.}$$
Moreover, if $v_1, \cdots , v_m$ is a basis of $V^r$ then $Z= Z^{r-1}$, $Z$ has a basis $z_1, \cdots , z_m$ with $dz_i\in A_{PL}(Y)$, and $[dz_i] = H(f)v_i$.

\vspace{3mm}\noindent {\bf Example 2: Nilpotent Serre fibrations}

Suppose $Y$ is a path connected CW complex and $Y = X_0 \stackrel{p_1}{\leftarrow} X_1\leftarrow \cdots$ is a possibly infinite sequence of elementary principal fibrations defined by homotopy pullback diagrams
$$\xymatrix{X_r \ar[rr] \ar[d]^{p_r} && P\ar[d]\\
X_{r-1}\ar[rr] && K(G_{r+1}, k_{r+1})\,, & r\geq 1}\,.$$
The resulting map $p : X:= \varprojlim_r X_r \to Y$ is a Serre fibration. If in addition, for each $n$ there is an $ r(n)$ such that $k_r>n$ if $r\geq r(n)$, then $p$ is a \emph{nilpotent fibration}. In this case the fibre of $X\to X_{r(n)}$ is $n$-connected and it follows that $H(X)= \varinjlim_r H(X_r)$. Moreover, the fibre, $F$, of $p$ is path connected and nilpotent, and $\pi_1(Y)$ acts nilpotently in $H(F)$. 

Now suppose in such a nilpotent fibration that each $G_r\otimes \mathbb Q$ is finite dimensional. Then there are quasi-isomorphisms
$$(A_{PL}(X_{r-1})\otimes \land Z_r,d) \stackrel{\simeq}{\longrightarrow } A_{PL}(X_r)$$
in which $d : Z_r\to A_{PL}(X_{r-1})$, dim$\, Z_r=$ dim$\, G_{r+1}\otimes \mathbb Q$, and $Z_r=Z_r^{k_r}$. These  combine to give a commutative diagram
$$\xymatrix{
(A_{PL}(Y)\otimes \land Z,d) \ar[rr]^\simeq && A_{PL}(X)\\
\mbox{} & A_{PL}(Y) \ar[lu] \ar[ru]_{A_{PL}(p)}}$$
in which  $(A_{PL}(Y)\otimes \land Z,d)$ is a   Sullivan extension. Here $Z = \oplus_r Z_r$ is a graded vector space of finite type.

In this case the construction in \S 12 gives a map of nilpotent Serre fibrations
$$
\xymatrix{ 
\mbox{} & X\ar[ld] \ar[dd]^{h_X} && F\ar[ll]\ar[dd]^{h_F}\\
Y\\
\mbox{} & X'\ar[lu] && \vert\land Z, \overline{d}\vert \ar[ll] 
}$$
in which $h_X$ and $h_F$ induce isomorphisms of rational cohomology. In particular, $H(F)$ is a graded vector space of finite type and $h_F$ is a rationalization of $F$.

\vspace{3mm} \noindent {\bf Example 3.} Let $G$ be a compact connected Lie group and $H$ a closed connected subgroup. The homogeneous space $G/H$ is the the quotient of $G$ by right multiplication by $H$.   Let now  $K$ be another closed connected subgroup of $G$ acting on $G$ on the left. If this action induces a free action on $G/H$, the quotient is called a biquotient and is denoted by $K\backslash G/H$. There is then a homotopy pullback diagram (\cite{Sing}) 
$$\xymatrix{K\backslash G/H \ar[rr] \ar[d] && BH \ar[d]^{Bf}\\
BK \ar[rr]^{Bg} && BG,}$$
where $f$ and $g$ denote the natural injections of $H$ and $K$ in $G$. 

 Then with the notation  of Example 2 (p.\pageref{hospace}), a Sullivan model for $K\backslash G/H$ is given by $$(\land s^{-1}V_K\otimes \land s^{-1}V_H\otimes \land V_G,d)$$
where $dx= H^*(Bf)(x) - H^*(Bg)(x)$, $x\in V_G$. 

\vspace{3mm}\noindent {\bf Example 4.} Suppose $G$ is a compact connected Lie group and $p : P\to M$ is a smooth principal fibre bundle. This fits in a homotopy pullback diagram
$$\xymatrix{P\ar[rr]\ar[d]^p && EG\ar[d]\\M\ar[rr]^f && BG}$$
where $EG$ is contractible. In the notation of Example 2 (p.\pageref{hospace}), write $H(G) = \land V_G$ and $H(BG)= \land s^{-1}V_G$; let $s^{-1}a_1, \cdots , s^{-1}a_n$ be a basis of $s^{-1}V_G$. Denoting the tangent bundles by $T_M$ and $T_P$, and the pullback of $T_M$ to $P$ by $p^*T_M$, Chern-Weil theory (\cite{Cartan}) uses an appropriate embedding of $p^*T_M$ in $T_P$ to construct closed differential forms $\omega_1, \cdots , \omega_n$ on $M$ representing $H(f)a_1, \cdots , H(f)(a_n)$. It also provides a quasi-isomorphism
$$(A_{DR}(M)\otimes \land (a_1, \cdots , a_n), 0)\stackrel{\simeq}{\to} A_{DR)}(M)\,, \hspace{1cm} Da_i = \omega_i\,,$$
which is one of the earliest relative Sullivan models.

\vspace{3mm}\noindent {\bf Example 5.} The free loop space $LX=\mbox{Top}(S^1, X)$ is the pullback in the diagram
$$\xymatrix{LX \ar[rr] \ar[d] && X^{[0,1]}\ar[d]^{p_0,p_1}\\
X \ar[rr]^{\Delta} && X\times X}$$
where $\Delta$ is the diagonal map and $p_i : X^{[0,1]}\to X$ denotes the evaluation at the point $i$.
Using this construction in perhaps the earliest application of Sullivan models, Vigu\'e-Poirrier and Sullivan (\cite{VS}) obtain the model of $LX$.

\begin{theoremb} Let $(\land V,d)$ be the minimal Sullivan model for a simply connected space $X$, then the minimal Sullivan model for $LX$ is
$$(\land V\otimes \land sV,d)\,, $$ where $D(sv)= -sdv$. Here the suspension $s$ has been extended  to the derivation to $\land V\otimes \land sV$ which is zero in  $sV$.\end{theoremb}

Using this model they prove that the free loop space of a simply connected compact manifold $M$ whose cohomology algebra  $H^*(M;\mathbb Q)$ requires at least two generators has unbounded Betti numbers. By a theorem of Gromoll and Meier (\cite{GM}) this implies with respect to any Riemannian metric that $M$ has infinitely many geometrically distinct closed geodesics (\cite{VS}). 

Now let $n_T$ denote the number of geometrically distinct closed geodesics of length $\leq T$. In \cite{Gr}, Gromov proves that for a simply connected compact manifold with a generic metric  there are constants $a>0$ and $c>0$ such that
$$n_T \geq a\cdot \mbox{max}_{p\leq cT} \mbox{dim}\, H^p(LM;\mathbb Q)\,.$$
Thus exponential growth of the Betti numbers of $LM$ would imply exponential growth of $n_T$.

\vspace{2mm}\noindent {\bf Example 6.} Generalizing Example 5, let $N\subset X\times X$ be a path connected subspace, where $X$ is a simply connected CW complex, and consider the pullback diagram, $$\xymatrix{X_N \ar[d]\ar[rr] && X^{[0,1]}\ar[d]^{p_0, p_1}\\N\ar[rr]^{i_N} && X\times X\,,}$$
where $i_N$ is the inclusion. If $(\land W,d)$ and $(\land V,d)$ are respectively the minimal Sullivan models of $N$ and of $X$ then (\cite{GHVi}) $X^N$ has a (not necessarily minimal) Sullivan model of the form $(\land W\otimes \land sV,D)$.

When $N$ is the graph of a periodic map $g : X\to X$ then the Sullivan-Vigu\'e theorem generalizes: if $H(X)$ has at least two generators fixed by $g$ then the Betti numbers of $X_N$ are unbounded. This in turn implies that, if $g$ is an isometry of a rationally hyperbolic Riemannian manifold, then there are infinitely many geometrically distinct $g$-invariants geodesics. 

\vspace{3mm} There remain open problems dealing with fibre squares.
\begin{enumerate}
\item[1.] If $X$ is a simply connected, rationally hyperbolic space, do the Betti numbers of $LX$ grow exponentially ?
\item[2.] Is the image of the Hurewicz homomorphism for $LX$ always finite dimensional if $X$ is a simply connected finite CW complex ?
\item[3.] If the base of a principal bundle is formal, what condition on the classes $[\omega_i]$ make $P$ formal as well ?
\end{enumerate}

\section{Mapping spaces}

Suppose $Y$ is a finite connected CW complex and $X$ is a nilpotent path connected space. Denote by $(X^Y,f)$ the path component of a map $f : X\to Y$ in the space $X^Y$ of all maps $Y\to X$. Then   (\cite{Hil}),  $(X^Y,f)$   is a nilpotent space and $(X^Y,f)_{\mathbb Q} = ((X_{\mathbb Q})^Y,f)$ 

More generally, suppose $p: X\to Y$ is a nilpotent fibration. Then (\cite{Mol}) each path component $(\Gamma, \sigma)$ in the space $\Gamma$ of all cross-sections of $p$ is nilpotent, and its rationalization is a path component in the space of cross-sections of $X'\to Y$, where $X'$  is the "rationalization  along the fibres" of $X$ (\cite{May}), described in Example 2 of \S 14 when the fibre $F$ of $p$ has rational homology of finite type.

Finally, Haefliger (\cite{Hae}), following an idea of Sullivan   and relying on earlier work of Thom (\cite{Th}) constructs a minimal Sullivan model for a  path component $(\Gamma, \sigma)$ of the space of cross-sections of a nilpotent fibration $p: Y\to X$ under the following hypotheses:
\begin{enumerate}
\item[$\bullet$] $Y$ is a finite, nilpotent connected CW complex.
\item[$\bullet$] The fibre, $F$, is path connected and $H(F)$ has finite type.
\end{enumerate}

Explicitly, since $Y$ is finite and nilpotent its minimal Sullivan model $(\land W,d)$ is a graded vector space of finite type, and for simplicity we write $(\land W,d)= A$. Then, since the fibration satisfies the hypothesis of Theorem \ref{12.1}, it has a Sullivan model of the form $(A\otimes \land V,d)$ where $(\land V, \overline{d})$ is a minimal Sullivan model for $F$. In this tensor product decomposition the choice of subspace $V$ can be and is made so that a Sullivan representative of $\sigma$ vanishes on $V$. 

Now let $B = A^\#$ be the dual differential graded coalgebra, and defines an algebra morphism
$$\varphi : A\otimes \land V\to A\widehat{\otimes}\land (B\otimes V)$$
by $\varphi (a) = a\widehat{\otimes} 1$ and $\varphi (v) = \prod_i a_i\widehat{\otimes} (b_i\otimes v)$, $v\in V$, where $a_i, b_j$ are dual basis of $A$ and $B$. Then a cdga $(\land (B\otimes V),d)$ is uniquely determined by the condition $\varphi\circ d= d\circ \varphi$. Finally, divide $\land (B\otimes V)$ by the ideal generated by $(B\otimes V)_{\geq 0}$ and $d(B\otimes V)_0$,   noting that $d(B\otimes V)_0\subset ((B\otimes V)_{\leq 0}\,\,{\scriptstyle \land} \land (B\otimes V)) \oplus (B\otimes V)_1$. Thus the quotient is   a Sullivan algebra $(\land Z,d)$  and Haefliger's theorem reads

\begin{theoremb} \label{thha}{\rm (\cite{Hae})}. The Sullivan algebra $(\land Z,d)$ is a Sullivan model for     $(\Gamma, \sigma)$. \end{theoremb}

\vspace{3mm}\noindent {\bf Remark.} An essential step in the proof of Theorem \ref{thha} is that (with the notation of Example 2 in the previous section)
$$\xymatrix{X^Y_{r+1} \ar[rr] \ar[d]^{p_{r+1}^Y} && P^Y\ar[d]\\
X^Y_r \ar[rr] && K(G_{r+2}, k_{r+1})^Y}$$
is a homotopy pullback diagram, and that the fibre of $p_{r+1}^Y$ is $K(G_{r+2}, k_{r+1})^Y$. Moreover, for any $K(G,m)$ the space $K(G,m)^Y$ has the weak homotopy type of the product $\prod_{i=0}^m K(H^{m-i}(Y,G),i)$.

\vspace{3mm}\noindent {\bf Example.} Suppose $X$ is a finite simply connected CW complex. Then multiplication in $G$ makes $X^G$ into a principal $G$-bundle. Constructing a Sullivan model for the fibration $X^G\to X^G/G\to BG$ is an open question, solved by Burghelea and Vigu\'e-Poirrier in the case $G= S^1$ and $X^{S^1}$ is the free loop space (\cite{BV}).

\vspace{3mm}\noindent {\bf Remark.} If $f : Y\to X$ is a continuous map, then $(X^Y,f)$ is homotopy equivalent to the space of sections of the trivial fibration $p: X\times Y$ homotopic to the section $s$, defined by $s(y)= (y, f(y))$. On the other hand, if $s$ is a section for a  fibration $p: X\to Y$, then $(\Gamma, s)$ is the homotopy fibre of the map $(X^Y,s)\to (Y^Y, id)$ (\cite{KSS}).

\vspace{3mm} Now let $(\land V,d)$ be a minimal Sullivan algebra and $K$ be a nilpotent simplicial set with finitely many non-degenerate simplices. The space $(\vert \land V,d\vert)^{\vert K\vert}$ is the geometric realization of the simplicial set whose $n$-simplices are the continuous maps 
$K\times \Delta^n\to  \vert \land V,d\vert$. By the adjunction property between spatial realization and $A_{PL}$, this is the geometric realization of the simplicial set whose $n$ simplices are the morphisms of cdga's,  $(\land V,d)\to  A_{PL}(K)\otimes A_{PL}(\Delta^n)$ (\cite{Be1}).  As above, denote by $A$  the minimal Sullivan model of $K$. Then there is a sequence of equivalences of simplicial sets
$$\mbox{Cdga} ((\land V,d), A_{PL}(K)\otimes A_{PL}(\Delta^*)) \simeq \mbox{Cdga} ((\land V,d), A\otimes A_{PL}(\Delta^*))$$
$$\simeq \mbox{Cdga} ((\land (B\otimes V),d),   A_{PL}(\Delta^*)) = \vert (\land B\otimes V),d)\vert\,.$$
Proceed then as above to obtain a model for the path components of the mapping space (\cite{BrS}, \cite{BrS0}).

\section{Acyclic closures and loop spaces}

The inclusion of a base point in a connected CW complex $X$ converts to the path space fibration $PX\to X$ with fibre the based loop space $\Omega X$. In Moore's construction of $PX$, composition of loops, $\mu_X : \Omega X\times \Omega X\to \Omega X$, makes $\Omega X$ into an associative H-space acting on the right on $PX$. 

The analogue in Sullivan's theory is the \emph{acyclic closure}
$$\xymatrix{ (\land V\otimes \land U,d) \ar[rr]^\simeq_\varepsilon && \mathbb Q\\ (\land V,d) \ar[u]^{\lambda_V} \ar[rru]_{\varepsilon_V}}$$
of a minimal Sullivan algebra, $(\land V,d)$. Note that if $V^1 \neq 0$ then $U^0\neq 0$; the augmentation $\varepsilon$ is chosen so that $\varepsilon (U) = 0$.
In this minimal $\land$-extension there is a linear isomorphism $\alpha : U\stackrel{\cong}{\to} V$ of degree $1$, and the differential satisfies
$$d(1\otimes u) - \alpha (u)\otimes 1 \in V {\scriptstyle \land} \land^+(V\oplus U)\,, \hspace{1cm} u\in U\,.$$
 In particular, the quotient differential in $\land U$ is zero.  

Multiplication in $\Omega X$ and the inverse $f^{-1}(t) = f(-t)$   also have  analogues: a diagonal
$$\Delta : \land U\to \land U\otimes \land U$$
and an involution $\omega$ in $\land U$,  which make  $\land U$ into a commutative graded Hopf algebra. Moreover, the spatial realizations 
$$\vert\Delta\vert : \vert\land U\vert \times \vert \land U\vert \to \vert\land U\vert\,, \hspace{5mm} \vert \omega\vert : \vert\land U\vert \to \vert\land U\vert\,, \hspace{5mm} \vert\varepsilon\vert : *\to \vert\land U\vert$$
make $\vert\land U\vert$ into a topological group. Similar constructions identify 
$$\xymatrix{ \vert \land V,d\vert &\vert\land V\otimes \land U,d\vert\ar[l] & \vert\land U\vert\ar[l]}$$
as a principal $\vert\land U\vert$ fibre bundle.

Now suppose $\varphi_X : (\land V,d) \to A_{PL}(X)$ is a minimal Sullivan model, so that there is a commutative diagram
$$\xymatrix{
(\land V,d) \ar[rr]\ar[d]^{\varphi_X} && (\land V\otimes \land U,d) \ar[rr]\ar[d]^{\varphi_{PX}} && (\land U, 0)\ar[d]^{\varphi_{\Omega X}}\\
A_{PL}(X) \ar[rr] && A_{PL}(PX) \ar[rr] && A_{PL}(\Omega X)\,.}$$
Moreover, even through $H(\varphi_{\Omega X})$ may not be an isomorphism, $$A_{PL}(\mu_X)\circ \varphi_{\Omega X} \sim \varphi_{\Omega X}\otimes \varphi_{\Omega X}\circ \Delta\,.$$

\vspace{3mm}\noindent {\bf Remark.} The constructions above generalize to an action of $\Omega X$ on the fibre when a continuous map is converted to a fibration, and to a diagonal $(\land Z, \overline{d})\to (\land Z, \overline{d})\otimes (\land U,0)$ for any $\land$-extension $(\land V,d)\to (\land V\otimes \land Z,d)\to (\land Z, \overline{d})$.

 With the notation above, the diagonal $\Delta$ dualizes to a multiplication,
 $$(\land U)^\#\otimes (\land U)^\#\to (\land U)^\#\,,$$
 which makes $(\land U)^\#$ into a graded algebra. The natural morphism $m_{ U} : \land U \to A_{PL}\vert\land U\vert$ induces a morphism $\land U\to H(\vert\land U\vert)$ of Hopf algebras, which then dualizes to an inclusion of graded algebras
 $$H_*(\vert\land U\vert, \mathbb Q) \rightarrowtail (\land U)^\#\,.$$
 
 Now assume that $H^1(\land V,d)$ and each $V^p$, $p\geq 2$, are finite dimensional, and denote by $L$ the homotopy Lie algebra of $(\land V,d)$. Then the holonomy representation of $\widehat{UL_0}\otimes UL_{\geq 1}$ in $\land U$ yields the isomorphism of graded algebras,
 $$\eta_L : \widehat{UL_0}\otimes UL_{\geq 1} \stackrel{\cong}{\to} (\land U)^\#,$$
defined by $\eta_L(a)(\Phi) = \varepsilon (a\bullet \Phi)$. With $\pi_1$ denoting $\pi_1(\land V,d)$, $\eta_L$ restricts to an isomorphism
 $$\mathbb Q[\pi_1]\otimes UL_{\geq 1}\stackrel{\cong}{\longrightarrow} H_*(\vert \land U\vert;\mathbb Q)$$
 of graded Hopf algebras.

 \section{Depth}
 
 The \emph{depth of a connected CW complex}, $X$, is a new homotopy invariant, constructed via its minimal Sullivan model $(\land V,d)$ as follows. Denote by $L_V$ and by $L_\alpha$ the homotopy Lie algebras of $(\land V,d)$ and $(\land V_\alpha, d)$, where the $(\land V_\alpha,d)$ are the sub Sullivan algebras satisfying $V_\alpha \subset V$ and dim$\, V_\alpha <\infty$. Then the $L_\alpha$ form an inverse system and $L_V = \varprojlim_\alpha L_\alpha$. Similarly, the $\widehat{UL_\alpha}$ also form an inverse system of $L_\alpha$-modules, and so $\varprojlim_\alpha \widehat{UL_\alpha}$ is naturally an augmented $L_V$-module.
 
 Finally, there are natural maps $H_*(\Omega X;\mathbb Q)\to \widehat{UL_\alpha}$ which then give a morphism $H_*(\Omega X;\mathbb Q) \to \varprojlim_\alpha \widehat{UL_\alpha}$. This exhibits $\varprojlim_\alpha \widehat{UL_\alpha}$ as a completion of $H_*(\Omega X;\mathbb Q)$ and we denote 
 $$\widehat{H}(\Omega X) = \varprojlim_\alpha \widehat{UL_\alpha}\,.$$
 
 \vspace{3mm}\noindent {\bf Definition.} The \emph{depth} of a connected CW complex $X$ is the least $p$ (or $\infty$) such that
 $$\mbox{Ext}^p_{UL_V}(\mathbb Q, \widehat{H}(\Omega X))\neq 0\,.$$
 This definition is motivated by the original definition of the depth of an (graded or ungraded) algebra $A$ augmented to $\mathbb Q$, as the least $p$ (or $\infty$) such that
 $$\mbox{Ext}_A^p(\mathbb Q, A)\neq 0\,.$$

 On the other hand, since each $V_\alpha$ is finite dimensional the holonomy representation of $L_\alpha$ in the fibre $\land U_\alpha$ of the acyclic closure gives an isomorphism
 $$\widehat{UL_\alpha} \stackrel{\cong}{\longrightarrow} (\land U_\alpha)^\#$$ of $\widehat{UL_\alpha}$-modules. Passing to inverse limits defines an isomorphism
 $$\widehat{H}(\Omega X) \stackrel{\cong}{\longrightarrow} (\land U)^\#$$
 of $UL_V$-modules, where $\land U$ is the fibre of the acyclic closure of $(\land V,d)$, and the representation of $UL_V$ is the holonomy representation. In particular,
 $$\mbox{Ext}_{UL_V}(\mathbb Q, \widehat{H}(\Omega X)) = \left[ \mbox{Tor}^{UL_V}(\mathbb Q, \land U)\right]^\#,$$
 and so $$\mbox{depth}\, X = \mbox{least $p$ (or $\infty$) such that Tor}_p^{UL_V}(\mathbb Q, \land U)\neq 0\,.$$

\vspace{3mm} For connected CW complexes $X$ we have the fundamental relation (\cite{Malcev})
 $$\mbox{depth}\, X\leq \mbox{cat}\, X\,.$$
 Moreover (\cite{depth}) let rad$\, L_V$ be the \emph{radical} of $L_V$, i.e., the sum of its solvable ideals. If depth$\, X<\infty$ then rad$\, L_V$ is finite dimensional and 
 $$\mbox{dim(rad}\, L_V)_{even} \leq \mbox{depth}\, X\,.$$
 
 \vspace{3mm}\noindent {\bf Remarks 1.} If dim$\, H^1(\land V,d)$ and dim$\, V^p<\infty$, $p\geq 2$, then depth$\,X=$ least $p$ (or $\infty$) such that $\mbox{Ext}^p_{UL_V}(\mathbb Q, \widehat{UL_V})\neq 0$.
 
 \vspace{2mm}{\bf 2.} In particular, if $X$ is simply connected and $H(X)$ has finite type then
 $$\mbox{depth}\, X = \mbox{least $p$ (or $\infty$) such that Ext}^p_{UL_V}(\mathbb Q, UL_V)\neq 0;$$
 i.e., depth$\, X = $ depth$\, UL_V$. 
 
 \vspace{2mm}{\bf 3.} Suppose $L_V$ is the quotient of a free Lie algebra $\mathbb L(S)$ by relations $x_\alpha \in \mathbb L(S)^{r_\alpha}$. If $H(\land sL_V,\delta)$ is finite dimensional, then again
 $$\mbox{depth}\, X = \mbox{least $p$ (or $\infty$) such that Ext}^p_{UL_V}(\mathbb Q, UL_V)\neq 0\,. $$
 It is an open question   when in general depth$\, X=$ depth$\, UL_V$.

 \vspace{3mm}Next, for any discrete group $G$ there are two associated depths: the depth of the classifying space $K(G,1)$ and the depth of the group ring $\mathbb Q[G]$. It is an open question   when and how these two invariants are related, but in the two examples below they are in fact equal.
 
 \vspace{3mm}\noindent {\bf Example 1.} (\cite{Artin}) If $M_g$ is a compact orientable Riemann surface of genus $g\geq 2$ then $M_g= K(\pi_1(M_g), 1)$ and 
 $$\mbox{depth}\, M_g= 2 = \mbox{depth}\, \mathbb Q[\pi_1(M_g)]\,.$$
 
 \vspace{3mm}\noindent {\bf Example 2.} (\cite{Artin}) A right angled Artin group $G$ is the quotient of a free group on generators $g_1, \cdots , g_n$ by relations of the form $g_ig_j= g_jg_i$. For such groups it also holds that
 $$\mbox{depth}\, K(G,1) = \mbox{depth}\, \mathbb Q[G]\,.$$
 
\vspace{3mm}There are a number of open questions about the depth of a connected CW complex, $X$, with minimal Sullivan algebra $(\land V,d)$, acyclic closure $(\land V\otimes \land U,d)$ and homotopy Lie algebra $L$. They include
\begin{enumerate}
\item[1.] Does depth$\, X$ depends only on the dga equivalence class of $C^*(X;\mathbb Q)$ ?
\item[2.] What follows if depth$\, X=$ cat$\, (\land V,d)$ ? (When $X$ is $1$-connected with finite Betti numbers, this implies that depth$\, X=$ gl dim$\, UL$.) 
\item[3.] Are there two connected CW complexes $X$ and $Y$ with isomorphic homotopy Lie algebras but different depths ?
\item[4.]  If depth$\, X = 1$ and dim$\, L/[L,L]\geq 2$, must $L$ contain a free Lie algebra on two generators ? (This is known if $X$ is simply connected with finite Betti numbers.)
\item[5.] If dim$\, H^1(X)<\infty$, is there a more general relation between depth$\, X$ and depth$\, \mathbb Q[G_L]$. Alternatively, is there a more general relation between depth$\, X$ and the least $p$ such that Ext$^p_{\mathbb Q[G]}(\mathbb Q, \widehat{\mathbb Q[G_L]})\neq 0$ ?
\end{enumerate}

\vspace{3mm} Every finite type graded Lie algebra, $L$, is the rational homotopy Lie algebra of a simply connected space  with finite Betti numbers, $\vert C^*(L)\vert$. Two questions arise directly. First, with what conditions on $L$ can we find a space $X$ with cat$\, X<\infty$ and $L_X= L$.  A necessary condition is given by the depth, because depth$\, L\leq$ cat$\, X$. On the other hand, let $L$ be a complete graded Lie algebra such that $L/[L,L] $ is a  graded vector space of finite type. Then on what conditions does there exist a space $X$ with $L_X= L$ ?

\section{Configuration spaces}

Let $M$ denote a compact $m$-dimensional manifold. The space of ordered configurations of $k$ points in $M$ is the space
$$F(M,k) = \{ (x_1, \cdots , x_k) \in M^k\,\vert\, x_i\neq x_j \mbox{ for } i\neq j\,\}\,.$$
An example of Longoni and Salvatore (\cite{LSa}) shows that the homotopy type of $F(M,k)$ is not a homotopy invariant of $M$ when $M$ is not simply connected. It remains a problem whether for simply connected compact manifolds the homotopy type (resp. the rational homotopy type) of $F(M,k)$ depends only on the homotopy type (resp. the rational homotopy type) of $M$. 

Partial results have been obtained using PD models. A PD (Poincar\'e duality) model for a simply connected closed manifold is a cdga $(A,d)$ weakly equivalent to the minimal Sullivan model of $M$ and that is a Poincar\'e duality algebra: there is a linear isomorphism $\varepsilon : A^m\to \mathbb Q$ such that the morphism $\varphi : A^p\to (A^{n-p})^\#$ given by $\varphi (a)(b)= \varepsilon (ab)$ is an isomorphism. 
In \cite{LSt} Lambrechts and Stanley show that every simply-connected closed manifold admits a PD model. Note that by duality $A^0= \mathbb Q \cdot 1$.
The diagonal class $D_A$ is then a cycle in $A\otimes A$ defined by
$$D_A = \sum_{deg\, a_i} a_i\otimes a_i'$$
where $a_i$ and $a_i'$ are graded basis of $A$ with $a_i\cdot a_j'= \delta_{ij}\omega$, $\omega$ denoting the fundamental class; i.e, a element with $\varepsilon (\omega)= 1$.

Let $(A,d)$ be a PD model for a manifold $M$. Denote by $p_i : A\to A^k$ and $p_{ij} : A^2\to A^k$ the injections,
$$p_i(a) = 1\otimes \cdots \otimes a\otimes 1\otimes \cdots \otimes 1\,,\hspace{5mm} p_{ij}(a,b)= 1\otimes \cdots \otimes a\otimes \cdots \otimes b\otimes \cdots \otimes 1$$
injecting $a$ (resp. $a$  and $b$) in position $i$ (resp. in positions $i$ and $j$). 
Form   the cdga
$$F(A,k) = (A^{\otimes k} \otimes \land (x_{ij}, i,j= 1\cdots ,k)/I,d)$$
where $dx_{ij}= p_{ij}(D_A)$, $deg\, x_{ij}= m-1$ and where the ideal $I$ is generated by the relations
$$\left\{\begin{array}{l} x_{ij}-(-1)^mx_{ji}, x_{ij}^2\\
p_i(x)\cdot x_{ij}- p_j(x)\cdot x_{ij}\\
x_{ij}x_{ik}+ x_{jk}x_{ki}+ x_{ki}x_{ij}\end{array}\right.$$

A main conjecture asserts  that $F(A,k)$ is weakly equivalent to the minimal Sullivan model of $F(M,k)$. This conjecture has been solved by Lambrechts and Stanley for $k=2$ when $M$ is 2-connected \cite{LSconf}  and  independently, by Kriz \cite{Kr} and Totaro \cite{To}, with $A= H(M)$ for   simply  connected complex projective manifolds. Recently Campos and Willwacher,  and independently Idrissi, have proved that $F(A,k)\otimes \mathbb R$ is weakly equivalent to $A_{PL}(M)\otimes \mathbb R$ for 1-connected closed manifolds (\cite{CW}, \cite{I}).

On the other hand, for manifolds of dimension $\geq 3$, Church proves a stability theorem for the action of the symmetric group $\Sigma_n$ on  $ H_*(F(M,n);\mathbb Q)$ (\cite{Ch}). There is a \emph{branching rule} $R_n^{n+1}$ that transforms an irreductible representation of $\Sigma_n$ into a representation of $\Sigma_{n+1}$. This branching rule is then extended to all finite dimensional representations by $R_n^{n+1}(V\oplus W)= R_n^{n+1}(V)\oplus R_n^{n+1}(W)$. The stability theorem of Church states that $H^i(F(M,n+1)) = R_n^{n+1} \left( H^i(F(M,n)\right)$ for $n\geq 2 i$. It is a natural question  if there is an analogous stability result for the rational homotopy groups of $F(M,n)$ or for their homotopy types. 

\section{The complement of arrangements}

An arrangement ${\mathcal A}$ of subspaces is a finite set of linear subspaces $X_i$ in some $\mathbb C^n$. Its complement is the space $M = \mathbb C^n \backslash \cup X_i$. Suppose first that the $X_i$ are hyperplanes   and denote by $\varphi_i$   linear forms with kernel $X_i$. Then Brieskorn shows that $M$ is a formal space, and that $H^*(M;\mathbb R)$ is the subalgebra of $A_{DR}(M)$ generated by the forms $\omega_i = \frac{1}{2\pi} d(\log \varphi_i)$ (\cite{Bri}). 

The arrangement of hyperplanes is of fiber-type if there is a tower $M= M_n \stackrel{p_n}{\to} M_{n-1}\to \cdots \stackrel{p_2}{\to} M_1= \mathbb C^*$ where each $M_k$ is the complement of an arrangement of hyperplanes in $\mathbb C^k$, and where $p_k$ is the restriction of a linear map $\mathbb C^k\to \mathbb C^{k-1}$ whose fiber   is a copy of $\mathbb C$ with finitely many points removed. In that case the minimal Sullivan model $(\land V,d)$ of $M$ satisfies $V = V^1$ (\cite{Fa}).

Consider now the more general case when the $X_i$ are general linear subspaces of $\mathbb C^n$. In \cite{FY},   Feichner and Yuzvinsky associate to the arrangement a complex    $(D,d)$ where $D$ is the $\mathbb Q$-vector space with basis the subsets  $\sigma = \{X_1, \cdots , X_p\} $ of ${\mathcal A}$. Order   the elements of the arrangement so  that in the writing of $\sigma$, $X_1<\cdots <X_p$. Write $\vert\sigma\vert = p$ and $\cap \sigma= \cap_i X_i$. Then the differential $d$ is defined by
$$d(\sigma) = \sum_{\{ i\,\, \vert \cap \sigma =   \cap \{\sigma\backslash \{i\}\}} (-1)^i \sigma\backslash \{X_i\}\,.$$
Moreover a product on $(D,d)$ is defined by 
$$\sigma\cdot \tau = \left\{\begin{array}{ll}
 (-1)^{\varepsilon (\sigma, \tau)} \sigma \cup\tau & \mbox{if }codim \cap \sigma + \mbox{codim }\cap \tau = \mbox{codim }\cap (\sigma \cup \tau)\\
0 & \mbox{otherwise}\end{array}\right.$$
The cdga $(D,d)$ is then weakly equivalent to the minimal Sullivan model of the complement in $\mathbb C^n$ of the union of the subspaces of the arrangement ${\mathcal A}$ \cite{FY}.

The arrangement is called geometric if   the lattice formed by the intersections of the elements of ${\mathcal A}$ satisfies the following property: if $X\subset Y$ and there is no $Z$ with $X\subset Z\subset Y$, then for any $T$ in the lattice there is no $Z$ between $X\cap T$ and $Y\cap T$. A geometric lattice is rationally a formal space (\cite{FY}).

\section{The group of homotopy self-equivalences of a space}

Let $X$ be a simply connected finite CW complex of dimension $m$. Denote by $X_{\mathbb Q}$ its rationalization and by ${\mathcal E}(X) $ the group of homotopy classes of self-homotopy equivalences of $X$, and by ${\mathcal E}_\#^m(X)$ the subgroup generated by self-equivalences inducing the identity in homotopy in degrees $\leq m$.   Dror and Zabrodsky prove that ${\mathcal E}_\#^m(X)$ is a nilpotent group (\cite{DZ}) and Maruyama proves that the natural map ${\mathcal E}_\#^m(X) \to  {\mathcal E}_\#^m(X_{\mathbb Q})$ is a rationalization \cite{Ma}. Now ${\mathcal E}_\#^m(X_{\mathbb Q})$ is anti-isomorphic to the group ${\mathcal E}_\#^m(\land V,d)$ of homotopy classes of self-equivalences of the minimal Sullivan model of $X$ that induce the identity on $\pi_{\leq m}(\land V,d)$. It follows (\cite{FM}) that
$$\mbox{nil}\,  {\mathcal E}_\#^m (X_{\mathbb Q}) < \mbox{cat}\, X\,.$$ 

\vspace{3mm} On the other hand, Costoya and Viruel have proved that every finite group $G$ can be realized as ${\mathcal E}(X)$ for a rationally elliptic space $X$ (\cite{CV}). The result is based on the fact that every finite group is the group of automorphisms of a finite graph $G$, and that we can associate to every graph $G$ an elliptic space $X$ with ${\mathcal E}(X)= \mbox{Aut}\, G$. 

A space $X$ is called \emph{homotopically rigid} if ${\mathcal E}(X) = \{1\}$ and a 1-connected closed manifold  is  \emph{inflexible} if      the degrees of the continuous self-maps form a finite set. Rational homotopy enables the construction of such manifolds: up to certain dimension there exist  infinitely many 1-connected rationally elliptic inflexible compact smooth manifolds (\cite{Am}), and there exist  infinitely many non equivalent homotopically rigid spaces   \cite{CMV}.

\section{Topological complexity}

A \emph{motion planning algorithm} for a space $X$ is a continuous section of the fibration $ev:= (p_0, p_1) : X^{[0,1]}\to X\times X$. Such a section exists if and only if $X$ is contractible. A local section $s : V\to X^{[0,1]}$ is called a motion planning rule on $V$. The topological complexity $TC(X)$ is 1 less than the minimum number of open sets $V$ covering $X\times X$ which admit motion planning rules. This is a homotopy invariant and
$$\mbox{cat}\, X\leq TC(X)\leq 2\mbox{ cat}\, X\,.$$

When $X$ is a 1-connected formal space with finite Betti numbers, then (\cite{LM}) $TC(X_{\mathbb Q})$ is the cup length $c_H$ of the kernel of the multiplication law $H(X)\otimes H(X)\to H(X)$. For general spaces $c_H$ is a lower bound for $TC(X_{\mathbb Q})$. There is an algebraic characterization in terms of the minimal model $(\land V,d)$ of $X$ that resembles the definition of the category of a minimal Sullivan model  (\cite{Car}): let $\mu : (\land V,d)\times (\land V,d) \to (\land V,d)$ be the multiplication map. Then TC$(X_{\mathbb Q})$ is the least integer $n$ such that the projection $(\land V,d)\otimes (\land V,d)\to (\land V,d)^{\otimes 2}/ (\mbox{ker}\, \mu)^{n+1})$ admits a homotopy retraction.

An interesting approximation to topological complexity has been obtained by Grant, Lupton and Oprea (\cite{GLO}): Let $f: Y\to X$ and $g: Z\to X$ be continuous maps between simply connected finite CW complexes. Suppose $\pi_*(f)\otimes \mathbb Q$ and $\pi_*(g)\otimes \mathbb Q$ are injective with disjoint images, then cat$\, Y_{\mathbb Q} + \mbox{cat}\, Z_{\mathbb Q}\leq  \mbox{TC}\,  (X_{\mathbb Q})$. In particular, 
$$\mbox{cat} (X_{\mathbb Q})+ \mbox{cat} (Y_{\mathbb Q}) \leq \mbox{TC}\,((X\vee Y)_{\mathbb Q})\,.$$

\vspace{1cm} Rational homotopy has led to important progress in many other domains in topology and geometry:   the topology of torus actions (\cite{AP}), the topology of complex manifolds (\cite{Mo}, \cite{Ci}), the intersection cohomology for manifolds with singularities (\cite{CST}), string topology, and the topology of the space of embeddings of a manifold into another (work of Lambrechts, Turchin, Volic and others, see for instance \cite{LTV}).

 \section{DG Lie models}
 
 Sullivan's approach to rational homotopy theory has a counterpart in Lie algebras. By the work of Quillen (\cite{Q}) each simply connected space   $X$ has a Lie model ${\mathcal L}_X$ that is a differential graded Lie algebra $(\mathbb L (V), d)$ with $V = s^{-1}H_*(X;\mathbb Q)$. Moreover (\cite{Maj}, \cite{BaL}) the cochain algebra $C^*({\mathcal L}_X)$ is quasi-isomorphic to the minimal Sullivan model,  $(\land V,d)$, of $X$, and   $H_*({\mathcal L}_X) = L_V$.
 Similarly each continuous map $f : X\to Y$ between simply connected   CW complexes has a Lie model $\varphi : {\mathcal L}_X\to {\mathcal L}_Y $, and each operation that can be done with Sullivan models of simply connected spaces has a counterpart with Lie models (\cite{Ta}, \cite{Q}, \cite{BaL}).
 
 In particular, the homotopy Lie algebra of the mapping space $Map^*_f(X,Y)$ of pointed maps homotopic to $f$ is given by the homology of the space $\mbox{Der}_{\varphi} ({\mathcal L}_X, {\mathcal L}_Y)$ of $\varphi$-derivations of ${\mathcal L}_X$ into ${\mathcal L}_Y$ (\cite{LS3}),
 $$H_*\mbox{Der}_\varphi ({\mathcal L}_X, {\mathcal L}_Y)) \cong L_{\mbox{Map}_f^*(X,Y)}\,.$$
 
 This isomorphism has been interpreted in terms of Hochschild and Andr\'e-Quillen homology (\cite{BL}). For that we need to recall that when $A$ is  a differential graded algebra   and $M$ and $N$ are differential graded $A$-modules, then 
 $$\mbox{Ext}^p_A(M,N) := H^p(\mbox{Hom}_A(P,N))\,,$$
 where $P\to M$ is a semifree resolution. 
 A particular example is given by the Hochschild cohomology of a differential graded algebra $A$ with coefficients in a bimodule $M$, $$HH^*(A,M):=\mbox{Ext}_{A\otimes A^{op}}(A,M)$$ where $A^{op}$ is the opposite algebra $a\cdot_{op}b = (-1)^{deg\, a\cdot deg\, b} b\cdot a$ and the action of $A\otimes A^{op}$ on $A$ is given by 
  $$(a\otimes b)\cdot c= (-1)^{deg\, b\cdot deg\, c} a\cdot c\cdot b\,.$$
 
 The morphism $\varphi : {\mathcal L}_X\to {\mathcal L}_Y$ equips $U{\mathcal L}_Y$ with a structure of $U{\mathcal L}_X$-bimodule, and there is a canonical isomorphism
 $$HH^*(U{\mathcal L}_X, U{\mathcal L}_Y) \cong \mbox{Ext}_{U{\mathcal L}_X} (\mathbb Q, U{\mathcal L}_Y)\,.$$
The action of ${\mathcal L}_X$ on the right is the adjoint action. For this action $U{\mathcal L}_Y$ decomposes into a direct sum 
$$U{\mathcal L}_Y = \mathbb Q \oplus {\mathcal L}_Y \oplus (U{\mathcal L}_Y)_{(2)} \oplus \cdots $$
Therefore $\mbox{Ext}_{U{\mathcal L}_X}(\mathbb Q, {\mathcal L}_Y)$ is a direct factor of $\mbox{Ext}_{U{\mathcal L}_X}(\mathbb Q, U{\mathcal L}_Y)$. 
 
Now recall   that, since ${\mathcal L}_X$ is the free Lie algebra $\mathbb L(V)$,  $\mathbb Q$ has a nice presentation as a  $(T(V),d)= U{\mathcal L}_X$-module.
 There is a natural quasi-isomorphism $$\psi : (T(V)\otimes (\mathbb Q \oplus sV),d)\stackrel{\simeq}{\longrightarrow}  (\mathbb Q, 0)\,,$$
where $d(sv) = v - Sdv$ and $S$ is the linear map $T^+V\to TV\otimes sV$ defined by  $S(v_1\cdots v_r) = (-1)^{\sum_{i<r} deg\, v_i} v_1\cdots v_{r-1} \otimes sv_r$. 
Extending each linear map $sV\to {\mathcal L}_Y$ into a $\varphi$-derivation induces    isomorphisms (\cite{BL}, \cite{Gat})
 $$\mbox{Ext}_{U{\mathcal L}_X}(\mathbb Q, {\mathcal L}_Y) \cong H_*(\mbox{Der}_\varphi ({\mathcal L}_X, {\mathcal L}_Y)) \cong L_{\mbox{Map}^*_f(X,Y)}\,.$$
  
  \vspace{3mm} Now let $\varphi : C^*(L)\to (A,d)$ be a model of a pointed map $f : X\to Y$, with dim$\, A<\infty$, and let $\tau\in A\otimes L$ be the corresponding MC element (Theorem \ref{Geth}). Then a Lie model for Map$_*(X,Y;f)$ is given by the dgl 
  $$(A\otimes L, d_\tau)\,, \hspace{1cm} d_\tau (a\otimes \ell)= d(a)\otimes \ell  + (-1)^{deg a} a\otimes d\ell + [\tau, a\otimes \ell]\,.$$
  
\vspace{3mm} In many contexts it is interesting to determine the Lie model of a space. For instance it is an open problem  to determine the Lie model of the space  of configurations of points in a manifold.

\vspace{3mm} Finally, let $(\land V,d)$ be a general minimal Sullivan algebra and $(\land V\otimes \land U,d)$ its acyclic completion. By hypothesis $V$ and $U$ are quipped with filtrations $V(p)$ and $U(p)$ such that $dV(p)\subset V(p-1)$ and $dU(p) \subset \land (V(p)\oplus U(p-1))$. We denote by ${\mathcal L}$ the Lie algebra of derivations $\theta$ of $\land V\otimes \land U$ such that $\theta (V)= 0$ and for each $p$, $\theta (U(p) ) \subset V\otimes \land U(p-1)$. Then (\cite{Malcev}) ${\mathcal L}
$ is a dgl whose homology is $L_V$. This dgl can be understood as a Lie model corresponding to $(\land V,d)$, because when $V= V^{\geq 2}$ is a   graded vector space of finite type then $\mathcal L$ is quasi-isomorphic to the usual Lie model associated to $(\land V,d)$ (\cite{Der}).
 
Question: Is it possible to make the association $(\land V,d)\mapsto \mathcal L$ the basis of a correspondence between cdga's and dgl's. In particular, is it possible to reconstruct directly $(\land V,d)$ from $\mathcal L$.

   \vspace{1cm} Institut de Recherche en Math\'ematique et en Physique, Universit\'e Catholique de Louvain, 2, Chemin du CYclotron, 1348 Louvain-La-Neuve, Belgium, yves.felix@uclouvain.be
   
   \vspace{3mm} Department of Mathematics, Mathematics Building, University of Maryland, College Park, MD 20742, United States, shalper@umd.edu
   
\end{document}